\documentclass[12pt]{article}

\usepackage{amsmath}
\usepackage{amsfonts}
\usepackage{latexsym}
\usepackage{graphicx}

\newtheorem{proposition}{Proposition}[section]
\newtheorem{theorem}{Theorem}[section]
\newtheorem{corollary}{Corollary}[section]
\newtheorem{remark}{Remark}[section]

\date{ }

\begin{document}

\title{Periodic Hamiltonian systems in shape optimization problems
with Neumann boundary conditions}

\author{Cornel Marius Murea$^1$, Dan Tiba$^2$\\
{\normalsize $^1$ D\'epartement de Math\'ematiques, IRIMAS,}\\
{\normalsize Universit\'e de Haute Alsace, France,}\\
{\normalsize cornel.murea@uha.fr}\\
{\normalsize $^2$ Institute of Mathematics (Romanian Academy) and}\\ 
{\normalsize Academy of Romanian Scientists, Bucharest, Romania,}\\ 
{\normalsize dan.tiba@imar.ro}
}

\maketitle

\begin{abstract}
The recent approach based on Hamiltonian systems and the
implicit parametri\-za\-tion theorem, provides a general fixed domain
approximation method in shape optimization problems,
using optimal control theory. In previous works,
we have examined Dirichlet boundary conditions with
distributed or boundary observation. Here, we discuss
the case of Neumann boundary conditions, with a
combined cost functional, including both distributed
and boundary observation. Extensions to nonlinear
state systems are possible. This new technique allows simultaneous boundary and topological variations and we also report numerical experiments confirming the theoretical results.

\vspace{3mm}
\textbf{Key Words:} Hamiltonian systems, implicit
parametrizations, shape optimization, optimal control,
Neumann boundary conditions, boundary and topological variations

\vspace{2mm}
\textbf{MCS 2020:} 49J50; 49M20; 49Q10
\end{abstract}

\section{Introduction}
\setcounter{equation}{0}

Shape optimization has started its development  especially in the last quarter of the previous century and we just quote several
monographs devoted to this subject Pironneau
\cite{Pironneau1984}, 
Haslinger and Neittaanm\"aki \cite{Haslinger1996},
Sokolowski and Zolesio \cite{Sokolowski1992},
Delfour and Zolesio \cite{Delfour2001},
Neittaanm\"aki, Sprekels and Tiba \cite{NS_Tiba2006},
Bucur and Buttazzo \cite{Bucur2005},
Henrot and Pierre \cite{Henrot2005}, Allaire \cite{A}, where more details on the history of the subject and comprehensive references can be found.
It is to be noted that, in general, just certain variants
of boundary variations are taken into account, while
topological variations of the unknown domains are
frequently not investigated.

A typical example of shape optimization problem, defined
on a given family of domains $\Omega \in \mathcal{O}$ (in general, it is assumed that $\Omega \subset D$, a
prescribed bounded domain), has
the following structure:
\begin{eqnarray}
\min_{\Omega\in\mathcal{O}} \int_\Lambda j\left(\mathbf{x}, y_\Omega(\mathbf{x})\right)d\mathbf{x},
\label{1.1}\\
A\, y_\Omega = f \hbox{ in }\Omega,
\label{1.2}\\
B\,y_\Omega = 0  \hbox{ on }\partial\Omega
\label{1.3}
\end{eqnarray}
where $\Lambda$ may be $\Omega$ or some fixed given
subdomain $E\subset \Omega$, or $\partial\Omega$; and
$B$ is some boundary operator expressing the boundary
condition, $A$ is some differential operator,
$f\in L^p(D)$, $p>2$ is given and $j(\cdot,\cdot)$ is a
a Carath\'eodory function.
More constraints on the
unknown domains $\Omega$, or on the state $y_\Omega$, more general cost
functionals may be taken into account. Regularity assumptions
on $\Omega\in\mathcal{O}$, on $j(\cdot,\cdot)$, other
hypotheses, will be imposed as necessity appears.

Many geometric optimization
problems arise in mechanics: minimize the thickness, the
volume, the stresses, etc., in a plate, a beam, a curved rod in
dimension three, an arch, a shell. Due to the
formulation of the mechanical models, the geometric
characteristics of the object (thickness, curvature)
enter as coefficients in the governing differential system.
Consequently, such geometric optimization problems take the
form of an optimal control problem in a given domain,
with the control acting in the coefficients. See \cite{bst}, \cite{alst},
\cite{NS_Tiba2006} Ch VI,
where detailed presentations, including numerical
examples, may be found.

In fact, general shape optimization problems
(\ref{1.1})-(\ref{1.3}) have a similar structure with optimal
control problems, the difference being that the minimization
parameter is the unknown geometry itself, $\Omega\in\mathcal{O}$.
It is a natural question to find a method that
reduces/approximates general optimal design problems to/via optimal
control theory, and some examples already appear in the
classical monograph of Pironneau \cite{Pironneau1984}.
In the case of Dirichlet boundary conditions several
approaches have been developed \cite{NP_Tiba2009}, \cite{N_Tiba2012}, \cite{MT2019},
\cite{MT2019a} allowing both shape
and topology optimization. Essential ingredients are
functional variations that combine both aspects and the
recent implicit parametrization method based on the
representation of the geometry via iterated Hamiltonian
systems \cite{Tiba2013},\cite{N_Tiba2015},\cite{Tiba2018},\cite{Tiba2018a}.
It turns out that this approach is very general and
we show here that it works in the case of Neumann boundary
conditions as well. This remains true for the Robin boundary conditions, nonlinear equations, 
etc., but we do not examine now such questions.
The methodology is of fixed domain type and it has important
advantages at the numerical level: it avoids remeshing  and recomputing the mass
matrix in
each iteration of the algorithm. Related ideas are also applicable in free boundary problems,
see \cite{HMT2016}, \cite{HMT2018}, optimization and control \cite{Tiba2020}.

Concerning topological variations, we underline that the well known level set method
\cite{OF}, \cite{OS}, \cite{A}, \cite{MAJ} is essentially different from our approach.
In our method, while we also use level functions, no Hamilton-Jacobi equation is needed
and simple ordinary differential Hamiltonian systems can handle the unknown geometry
and its variations. We work in dimension two, $D \subset \mathbb{R}^2$, since the
important periodicity argument is based on the Poincare-Bendixson theorem
\cite{Hirsch2014}, \cite{Pon}, and certain related developments.
This is a case of interest in shape optimization.

The paper is organized as follows. In the next Section, we
collect some preliminaries and we give the precise
formulation of the problem. Both distributed and boundary
observations are taken into account. In Section 3 we
introduce the fixed domain approximation process as an
optimal control problem, we prove a general approximation property under very weak
conditions and we also obtain some error estimates. As a corollary of the employed
methods, an existence result is proved as well. Section 4 is devoted to the
differentiability properties of our approach, that give the basis for numerical
algorithms of gradient type. A key technical development is the proof of the
differentiability of the period in Hamiltonian systems, with respect to functional
variations. Discretization and numerical examples are discussed in the last two Sections.

\section{Problem formulation and preliminaries}
\setcounter{equation}{0}

Let $\mathcal{O}$ be a given family of open, connected sets,
$\Omega \subset D$, not necessarily simply connected, where
$D\subset\mathbb{R}^2$ is a bounded domain and $\Omega$,
$D$ have both $\mathcal{C}^{1,1}$ boundaries.

In each $\Omega\in \mathcal{O}$, we consider the Neumann
boundary value problem
\begin{eqnarray}
-\Delta y_\Omega + y_\Omega= f \hbox{ in }\Omega,
\label{2.1}\\
\frac{\partial y_\Omega}{\partial n} = 0 \hbox{ on }
\partial\Omega,
\label{2.2}
\end{eqnarray}
where $f\in L^p(D)$, $p>2$ is given. It is known that
(\ref{2.1}), (\ref{2.2}) has a unique solution
$y_\Omega\in W^{2,p}(\Omega)$, more general elliptic operators
may be taken into account in (\ref{2.1}) or the regularity
conditions on the boundary may be relaxed, Grisvard
\cite{Grisvard1985}.
Here, it is important to work in $\mathbb{R}^2$ since
Poincar\'e-Bendixson type arguments are essential in the
proof of the global existence result for the Hamilton
system (\ref{2.9})-(\ref{2.11}) that are introduced in the sequel for the description
of the unknown geometries.
In fact, all the other arguments to be used in this work
are valid in arbitrary dimension, where iterated Hamiltonian
systems are necessary for the description of the geometry
and their solution is local \cite{Tiba2018}.

We associate to the system (\ref{2.1}), (\ref{2.2}) a cost
functional that combines distributed and boundary observation
(the necessary regularity conditions are detailed in the sequel):
\begin{equation}\label{2.3}
\min_{\Omega\in\mathcal{O}}
\left\{
\int_E J\left(\mathbf{x}, y_\Omega(\mathbf{x})\right)
d\mathbf{x}
+\int_{\partial\Omega}
j\left(\mathbf{x}, y_\Omega(\mathbf{x})\right)d\sigma
\right\},
\end{equation}
where $E\subset\subset D$ is a given subdomain such that
$E\subset\Omega$ for any $\Omega\in\mathcal{O}$ and 
$J(\cdot,\cdot)$, $j(\cdot,\cdot)$ are Carath\'eodory
functions. More restrictions (for instance, on the state
$y_\Omega$) may be added to the shape optimization problem
(\ref{2.1})-(\ref{2.3}), denoted by $(\mathcal{P})$.
More assumptions will be formulated as necessity appears.

The approach based on functional variations \cite{N_Tiba2012}, \cite{NP_Tiba2009},
\cite{Tiba2018a} assumes
that the family of admissible domain $\mathcal{O}$ is
obtained starting from a family
$\mathcal{F}\subset\mathcal{C}(\overline{D})$ of level
functions via the relation:
\begin{equation}\label{2.4}
  \Omega=\Omega_g=int\left\{ \mathbf{x}\in D;\ g(\mathbf{x})
  \leq 0\right\},\quad g\in \mathcal{F}.
\end{equation} 

While $\Omega_g$ defined in (\ref{2.4}) is an open set and
may have many connected components, the domain $\Omega_g$
that we use in the sequel is the component that contains
$E$. This is possible if we assume
\begin{equation}\label{2.5}
g(\mathbf{x})\leq 0,\quad \forall \mathbf{x}\in E,
\quad \forall g \in \mathcal{F}.
\end{equation}
Another variant, possible to be used in the definition of
the domain $\Omega_g$, is to assume that
\begin{equation}\label{2.6}
\mathbf{x}_0\in \partial\Omega_g,
\quad \forall g \in \mathcal{F}
\end{equation}
for some $\mathbf{x}_0\in D\setminus\overline{E}$, given. One has to impose on
the family $\mathcal{F}$ the simple constraint

\begin{equation}\label{2,5}
g(\mathbf{x}_0) = 0,
\quad \forall g \in \mathcal{F}.
\end{equation}

In this context, it is important to consider the closed
bounded set:
\begin{equation}\label{2.7}
G=\left\{  \mathbf{x}\in D;\ g(\mathbf{x})=0\right\}
\end{equation}
associated to any $g \in \mathcal{F}$.
If $\mathcal{F}\subset\mathcal{C}(\overline{D})$ without
further conditions, then $meas(G)>0$ is possible.
We further assume, see \cite{Tiba2018a}, that
$\mathcal{F}\subset\mathcal{C}^1(\overline{D})$ and
\begin{equation}\label{2.8}
  |\nabla g( \mathbf{x})| > 0,\quad \forall\mathbf{x}\in G,
  \quad\forall g \in \mathcal{F}.
\end{equation}
Then, by (\ref{2.6})-(\ref{2.8}) and the implicit functions
theorem, we get $G=\partial\Omega_g$ and the Hamiltonian
system
\begin{eqnarray}
z_1^\prime(t) & = & -\frac{\partial g}{\partial x_2}\left(z_1(t),z_2(t)\right),\quad t\in I_g,
\label{2.9}\\
z_2^\prime(t) & = &  \frac{\partial g}{\partial x_1}\left(z_1(t),z_2(t)\right),\quad t\in I_g,
\label{2.10}\\
\left(z_1(0),z_2(0)\right)& = & \mathbf{x}_0
\in \partial\Omega_g,
\label{2.11}
\end{eqnarray}
where $I_g$ is the local existence interval for
(\ref{2.9})-(\ref{2.11}), gives a local parametrization of
$\partial\Omega_g$ around $\mathbf{x}_0$, \cite{Tiba2013}.
The solution is unique due to the Hamiltonian structure \cite{Tiba2018}.
We also assume that
\begin{equation}\label{2.12}
  g( \mathbf{x}) > 0,\quad \forall\mathbf{x}\in \partial D,
  \quad\forall g \in \mathcal{F}
\end{equation}
which ensures that $G\cap \partial D=\emptyset$ for
$g \in \mathcal{F}$.

Notice that the family $\mathcal{O}$ of domains defined by
(\ref{2.4})-(\ref{2.5}) is very rich,
they may be multiply connected and this is one reason
why the above approach combines boundary and topological
variations in shape optimization.

Moreover, under hypothesis (\ref{2.8}), we get
$\partial\Omega_g$ of class $\mathcal{C}^1$ and more
regularity can be obtained if more regularity is imposed
on $\mathcal{F}$. This ensures the previously mentioned
regularity properties for the solution of 
(\ref{2.1}), (\ref{2.2}) and the cost (\ref{2.3}) and
its approximation (in the next section), are well defined.

It is proved in \cite{Tiba2018a}, that hypotheses (\ref{2.8}) and
(\ref{2.12}) are sufficient for the global existence
in (\ref{2.9})-(\ref{2.11}).

\begin{theorem}\label{theo:2.1}
For any $\mathbf{x}_0\in D\setminus E$, the solution of
(\ref{2.9})-(\ref{2.11}) is periodic and $I_g$ may
be chosen as its period, $I_g = [0,T_g]$.
\end{theorem}  

Namely, the limit cycle situation from the
Poincar\'e-Bendixson theory is not possible here. If
$\partial\Omega_g$ is not connected, its complete
description may be obtained via
(\ref{2.9})-(\ref{2.11}), by choosing an initial condition
on each component. Another crucial property proved in \cite{Tiba2018a} is

\begin{theorem}\label{theo:2.2}
Under the above hypotheses, the compact set $G$ has
a finite number of connected components, for any fixed
$g \in \mathcal{F}$.
\end{theorem} 

\noindent
Clearly, the number of the connected components may be
unbounded over the whole $\mathcal{F}$.

\section{Approximation and existence}
\setcounter{equation}{0}

The approximation of shape optimization problems via
cost penalization was  introduced in \cite{Tiba2018a} and further
developed in \cite{MT2019}. The idea is to penalize the
boundary condition on the unknown domains. This is
possible due to the Hamiltonian representation of the
unknown geometries, Thm. \ref{theo:2.1} and Thm.
\ref{theo:2.2}. We use here a penalization variant that
has good differentiability properties and is formulated as
an optimal control problem ($\epsilon>0$):
\begin{eqnarray}
&&\min_{g,u}
\left\{
\int_{E} J\left(\mathbf{x},y(\mathbf{x})\right)
d\mathbf{x}
+
\int_{I_g}
j\left(\mathbf{z}(t), y(\mathbf{z}(t))\right)
\sqrt{ (z_1^\prime(t))^2 + (z_2^\prime(t))^2} dt
\right.
\nonumber\\
&&
\left.
+\frac{1}{\epsilon}
\int_{I_g}
\left[
  \nabla y(z_1(t),z_2(t))\cdot
  \frac{\nabla g(z_1(t),z_2(t)) }{
  |\nabla g(z_1(t),z_2(t)) |}
\right]^2
\sqrt{ (z_1^\prime(t))^2 + (z_2^\prime(t))^2} dt
\right\}
\label{3.1}
\end{eqnarray}
subject to
\begin{eqnarray}
-\Delta y + y& = & f+ g_+^2 u,\quad\hbox{in } D,
\label{3.2}\\
y & = & 0,\quad\hbox{on } \partial D,
\label{3.3}
\end{eqnarray}
and (\ref{2.5}).
Above $\mathbf{z}(t)=(z_1(t),z_2(t))$ is the solution of
(\ref{2.9})-(\ref{2.11}), the state $y\in W^{2,p}(D)
\cap H_0^1(D)$ from (\ref{3.2}), (\ref{3.3}) clearly
depends on $g \in \mathcal{F}$ and $u$ is measurable
such that $g_+^2 u \in L^p(D)$, $p>2$.
In dimension 2, we have $y\in \mathcal{C}^1(\overline{D})$
by the Sobolev theorem and all the terms in (\ref{3.1})
make sense.
The penalization term in (\ref{3.1}) is a detailed
formula for
$$
\int_{\partial \Omega_g}
\left|
\frac{\partial y}{\partial n}
\right|^2 d\sigma
$$
based on the Hamiltonian representation
(\ref{2.9})-(\ref{2.11}) of $\partial \Omega_g$ and the fact
that the unit normal to $\partial \Omega_g=G$ is given by
$\frac{\nabla g(z_1(t),z_2(t)) }{
  |\nabla g(z_1(t),z_2(t)) |}$ in $(z_1(t),z_2(t))\in \partial
\Omega_g$ and it is well defined due to (\ref{2.8}). In case
$\partial\Omega_g$ has several connected components (their
number is finite by Thm. \ref{theo:2.2}) then the penalization
term is replaced by a finite sum of similar terms, with some
initial condition in (\ref{2.9})-(\ref{2.11}) fixed on each
component.
It is to be noticed that, in the ``extended'' equation
(\ref{3.2}), (\ref{3.3}), we have Dirichlet boundary conditions,
while the original state system (\ref{2.1}), (\ref{2.2}) is a
Neumann boundary value problem. It turns out that the
approximation properties of (\ref{3.1})-(\ref{3.3}) remain valid
even with this change of boundary conditions and we want to
stress this property. In fact, it is also easier to work
with (\ref{3.3}) in the finite element discretization, in the
next sections.

\begin{proposition}\label{prop:3.1}
Let $J(\cdot,\cdot)$ and $j(\cdot,\cdot)$ be 
Carath\'eodory functions on $D\times\mathbb{R}$, bounded from
below by a constant and let
$\mathcal{F}\subset \mathcal{C}^2(\overline{D})$ satisfy
(\ref{2.8}), (\ref{2.12}). Denote by $[y_n^\epsilon ,g_n^\epsilon ,u_n^\epsilon]$ a minimizing sequence in the penalized problem (\ref{3.1})-(\ref{3.3}), (\ref{2.5}). Then, on a subsequence
denoted by $n(m)$ the pairs
$[\Omega_{g_{n(m)}^\epsilon}, y_{n(m)}^\epsilon]$ (not necessarily
admissible) give a minimizing cost in (\ref{2.3}),
satisfy (\ref{2.1}) and (\ref{2.2}) is valid with a perturbation
of order $\epsilon^{1/2}$.
\end{proposition}

\noindent
\textbf{Proof.}
The proof follows the ideas from \cite{Tiba2018a}, \cite{MT2019}.
Let $[y_{g_m}, g_m]\in W^{2,p}(\Omega_{g_{m}})\times
\mathcal{F}$ be a minimizing sequence for the problem 
(\ref{2.1})-(\ref{2.5}). Here, $\partial \Omega_{g_{m}}$ is
$\mathcal{C}^2$ and this ensures the regularity
$y_{g_{m}}\in W^{2,p}(\Omega_{g_{m}})$ due to $f\in L^p(D)$.
There is $\widetilde{y}_{g_{m}}\in
W^{2,p}(D\setminus \overline{\Omega}_{g_{m}})$, not unique,
such that $\widetilde{y}_{g_{m}}=y_{g_m}$ on
$\partial \Omega_{g_{m}}$, 
$\frac{\partial \widetilde{y}_{g_{m}}}{\partial \mathbf{n}}
= \frac{\partial y_{g_{m}}}{\partial \mathbf{n}} = 0$
on $\partial \Omega_{g_{m}}$, $\widetilde{y}_{g_{m}}=0$ on
$\partial D$.
We define an admissible control in (\ref{3.2}) by 
\begin{equation}\label{3.4}
u_{g_{m}} =
-\frac{\Delta \widetilde{y}_{g_{m}} + f -\widetilde{y}_{g_{m}}}
{(g_m)_+^2},
\quad\hbox{in } D\setminus \overline{\Omega}_{g_{m}},
\end{equation}
and zero otherwise.
We infer by (\ref{3.4}) that $(g_m)_+^2 u_{g_m}$ is in
$L^p(D)$ and $g_m$, $u_{g_{m}}$ is an admissible
control pair for the penalized problem
(\ref{3.1})-(\ref{3.3}), (\ref{2.5}). Moreover, 
the corresponding state in (\ref{3.2}) is obtained by
concatenation of $y_{g_{m}}$ and $\widetilde{y}_{g_{m}}$ and
the corresponding penalization term in (\ref{3.1}) is null.
That is the corresponding costs in (\ref{3.1}) and in (\ref{2.3})
are the same. This construction is also valid in the case
$\Omega_{g_{m}}$ is not simply connected.

We obtain
\begin{eqnarray}
&&
\int_{E} J\left(\mathbf{x},y_{n(m)}^\epsilon(\mathbf{x})\right)
d\mathbf{x}
+
\int_{I_{g_{n(m)}}}
  j\left(\mathbf{z}_{n(m)}(t),
  y_{n(m)}^\epsilon(\mathbf{z}_{n(m)}(t))\right)
| \mathbf{z}_{n(m)}^\prime(t) | dt
\nonumber\\
&&
+\frac{1}{\epsilon}
\int_{I_{g_{n(m)}}}
\left[
  \nabla y_{n(m)}^\epsilon(\mathbf{z}_{n(m)}(t))\cdot
  \frac{\nabla g_{n(m)}^\epsilon(\mathbf{z}_{n(m)}(t)) }{
  |\nabla g_{n(m)}^\epsilon(\mathbf{z}_{n(m)}(t)) |}
\right]^2
| \mathbf{z}_{n(m)}^\prime(t) | dt
\nonumber\\
&\leq &
\int_{E} J\left(\mathbf{x},y_m(\mathbf{x})\right)
d\mathbf{x}
+\int_{\partial\Omega_{g_m}} j\left(\mathbf{x},y_m(\mathbf{x})\right)
d\sigma \rightarrow \inf (\mathcal{P})
\label{3.5}
\end{eqnarray}
for $m\rightarrow\infty$.
In (\ref{3.5}), the index $n(m)$ is big enough in order
to have the inequality valid and $\mathbf{z}_n$ is the
solution of (\ref{2.9})-(\ref{2.11}) associated to
$g_n^\epsilon$ (for simplicity, we don't write
$\mathbf{z}_n^\epsilon$).

Since $J\left(\cdot,\cdot\right)$ and
$j\left(\cdot,\cdot\right)$ are bounded from below by
constants, from (\ref{3.5}), we get the boundedness of the
penalization term on the subsequence $n(m)$. This yields
the last statement of Proposition \ref{prop:3.1}, on
$\partial\Omega_{g_{n(m)}^\epsilon}$.
As $\left(g_{n(m)}^\epsilon\right)_+$ is null in
$\Omega_{g_{n(m)}^\epsilon}$, we see that (\ref{2.1}) is
satisfied in $\Omega_{g_{n(m)}^\epsilon}$, due to (\ref{3.2}).
The minimizing property of the sequence
$\left[ \Omega_{g_{n(m)}^\epsilon}, y_{n(m)}^\epsilon\right]$ in
the original cost (\ref{2.3}) is again an obvious
consequence of (\ref{3.5}), by the positivity of the
penalization term(s).\quad$\Box$

\noindent
By the Weierstrass theorem, there is $m_g > 0$ such that (\ref{2.8}) becomes

\begin{equation}\label{3,6}
  |\nabla g( \mathbf{x})| \geq  m_g,\quad \forall\mathbf{x}\in G,
  \quad\forall g \in \mathcal{F}.
\end{equation}

\noindent
In order to strengthen the approximation property in Proposition \ref{prop:3.1},
we impose that $\mathcal{F} $ is bounded in $\mathcal{C}^2 (\overline{D})$ and we
require  uniformity in (\ref{2.8}), (\ref{3,6}), where $m > 0$ is some given constant: 
\begin{equation}\label{3,7}
  |\nabla g( \mathbf{x})| \geq  m,\quad \forall\mathbf{x}\in G,
  \quad\forall g \in \mathcal{F}.
\end{equation}

\noindent
Notice that (\ref{3,7}) or the boundedness of $\mathcal{F} $ don't modify
the topological characteristics of the family of admissible domains
$\Omega_g, \; g \in \mathcal{F}$.
We denote by $y_{n, \epsilon}$ the solution of (\ref{2.1}), (\ref{2.2})
in $\Omega_{g_n^\epsilon}$.

\begin{proposition}\label{prop:3,2}
Under the above assumptions, there is an absolute constant $C > 0$ such that
$$
|y_{n, \epsilon} - y_n^\epsilon |_{H^1 (\Omega_{g_n^\epsilon})} \leq C \epsilon^{1/4}.
$$
\end{proposition}

\noindent
\textbf{Proof.}
We take the difference of the equations (\ref{2.1}) in $\Omega_{g_n^\epsilon}$
corresponding to $y_{n, \epsilon}, \; y_n^\epsilon$ and we multiply by
$y_{n, \epsilon} - y_n^\epsilon$. Then, we get:

$$
|y_{n, \epsilon} - y_n^\epsilon |^2_{H^1 (\Omega_{g_n^\epsilon})} = -\int_{\partial \Omega_{g_n^\epsilon}}
(\frac{\partial  y_n^\epsilon}{\partial n})
 (y_{n, \epsilon} - y_n^\epsilon) d\sigma \leq c\epsilon^{1/2} 
 |y_{n, \epsilon} - y_n^\epsilon|_{L^2(\partial \Omega_{g_n^\epsilon})}, 
$$

\noindent
where $c > 0$ is an absolute constant corresponding to the evaluation
of the penalization term in (\ref{3.1}), from the last statement
in Proposition \ref{prop:3.1}.

By (\ref{3,7}) and Green's formula, we have:
\begin{eqnarray*}
&&
m|y_{n, \epsilon} - y_n^\epsilon|^2_{L^2(\partial \Omega_{g_n^\epsilon})}
\leq  \int_{\partial \Omega_{g_n^\epsilon}} |y_{n, \epsilon} - y_n^\epsilon|^2 \;
|\nabla g_n^\epsilon |d \sigma
=\int_{\partial \Omega_{g_n^\epsilon}} |y_{n, \epsilon} - y_n^\epsilon|^2
\nabla g_n^\epsilon \cdot \nu_\epsilon d \sigma \\
&&
\leq \int_{\Omega_{g_n^\epsilon}} |y_{n, \epsilon} - y_n^\epsilon|^2 |\Delta g_n^\epsilon|dx
+ 2\int_{\Omega_{g_n^\epsilon}} |y_{n, \epsilon} - y_n^\epsilon|
|\nabla(y_{n, \epsilon} - y_n^\epsilon) \cdot \nabla g_n^\epsilon|dx\\
&&
\leq M[|y_{n, \epsilon} - y_n^\epsilon|^2_{L^2(\Omega_{g_n^\epsilon})}
+ |y_{n, \epsilon} - y_n^\epsilon|_{L^2(\Omega_{g_n^\epsilon})}
|\nabla(y_{n, \epsilon} - y_n^\epsilon)|_{L^2(\Omega_{g_n^\epsilon})}]
\leq M[|y_{n, \epsilon} - y_n^\epsilon|^2_{L^2(\Omega_{g_n^\epsilon})} \\
&&  
+ \epsilon^{1/2} |\nabla(y_{n, \epsilon} - y_n^\epsilon)|^2_{L^2(\Omega_{g_n^\epsilon})}
+ \epsilon^{-1/2}|y_{n, \epsilon} - y_n^\epsilon|^2_{L^2(\Omega_{g_n^\epsilon})} ],
\end{eqnarray*}

\noindent
where we also use the binomial inequality (with the same $\epsilon$ as
in Proposition \ref{prop:3.1}) together with the boundedness of
$\mathcal{F}$ in $\mathcal{C}^2 (\overline{D})$. The notation
$\nu_\epsilon$ is the normal to the domain $\Omega_{g_n^\epsilon}$. 

\noindent
Combining the above two inequalities, we end the proof.\quad$\Box$

\begin{remark}\label{rem:3.1}
We note the very weak hypotheses on the cost functional in
Proposition \ref{prop:3.1}. Together with Proposition \ref{prop:3,2},
the justification for the use of the
control problem (\ref{3.1})-(\ref{3.3}), (\ref{2.5}) in
the approximation of $(\mathcal{P})$, is obtained.
A detailed study of the convergence properties when
$\epsilon \rightarrow 0$, for a distributed cost functional,
is performed in \cite{Tiba2018a}.
\end{remark}  

\begin{corollary}
Under assumption (\ref{3,7}) and the boundedness of $\mathcal{F}$ in
$\mathcal{C}^1 (\overline{D})$, the shape optimization problem has
at least one optimal solution $\Omega^*$.
\end{corollary}

\noindent
\textbf{Proof.}
Condition (\ref{3,7})  allows to apply the implicit function theorem
around any point $(x,y) \in G$ and to obtain the local representation
of $G$ via some function $y = y(x) $. In particular, also taking into
account the boundedness of $\mathcal{F}$ in $\mathcal{C}^1 (\overline{D})$,
it  yields that $ y^\prime (x) = - \frac{g_x(x,y(x))}{g_y (x,y(x))}$ is bounded,
uniformly with respect to the family of admissible domains, under
appropriate choices of the local axes. This allows the application
of well known existence results due to Chenais (see \cite{Pironneau1984},
Ch. 3.3)  and to end the proof.\quad$\Box$

\section{Directional derivative}
\setcounter{equation}{0}

We consider now functional variations $g+\lambda r$, $u+\lambda v$,
$r \in \mathcal{F}$, $\lambda \in \mathbb{R}$, $v\in L^p(D)$.
In the sequel, we shall take into account the condition
(\ref{2.6}), (\ref{2,5}) for $g$, $r$ in the identification of the corresponding domains 
from (\ref{2.4}). This is also necessary in
(\ref{2.9})-(\ref{2.11}) and at the numerical level it is
very easy to implement (finding some $\mathbf{x}_0$ arises to solve $g(\mathbf{x}) = 0$, which is a standard routine, and to use (\ref{2.9})-(\ref{2.11}) to identify such initial conditions on each connected component of $G$ by elimination; see \cite{MT2019} for other details). Notice that the perturbations of $u$ are always admissible since we have no constraints on $u$ and the perturbations of $g$ satisfy (\ref{2,5}), (\ref{2.8}), (\ref{2.12}) for $|\lambda|$ small enough (depending on $g$).

We denote by $y_\lambda\in W^{2,p}(D)$,
$\mathbf{z}_\lambda\in \mathcal{C}^1(\mathbb{R})$ the solutions
of (\ref{3.2}), (\ref{3.3}) and (\ref{2.9})-(\ref{2.11})
corresponding to the above variations, respectively. From
the previous section, we know that $\mathbf{z}_\lambda$ is
periodic with some period $T_\lambda>0$ and we take its
definition interval to be $[0,T_\lambda]$.  In \cite{MT2019}, it is
proved under conditions (\ref{2.8}), (\ref{2.12}), that
$T_\lambda\rightarrow T$ as $\lambda\rightarrow 0$, where $T$
is the period of $\mathbf{z}$, i.e. $I_g=[0,T]$.

\begin{proposition}\label{prop:3.2}
  The system in variations corresponding to (\ref{3.2}),
  (\ref{3.3}),
(\ref{2.9})-(\ref{2.11}) is:
\begin{eqnarray}
-\Delta q + q & = & g_+^2 v + 2g_+ u\,r,\quad\hbox{in } D,
\label{3.6}\\
q & = & 0,\quad\hbox{on } \partial D,
\label{3.7}\\
w_1^\prime& = & -\nabla\partial_2 g(\mathbf{z})\cdot \mathbf{w}
-\partial_2 r(\mathbf{z}),\quad\hbox{in } [0,T],
\label{3.8}\\
w_2^\prime& = & \nabla\partial_1 g(\mathbf{z})\cdot \mathbf{w}
+\partial_1 r(\mathbf{z}),\quad\hbox{in } [0,T],
\label{3.9}\\
w_1(0)& = &0,\ w_2(0) = 0,
\label{3.10}
\end{eqnarray}
where $q=\lim_{\lambda\rightarrow 0}\frac{y_\lambda-y}{\lambda}$, 
$\mathbf{w}=[w_1,w_2]=\lim_{\lambda\rightarrow 0}
\frac{\mathbf{z}_{\lambda} - \mathbf{z}}{\lambda}$
and the limits exists in $W^{2,p}(D)$, 
respectively $\mathcal{C}^1([0,T])$. 
\end{proposition}

\noindent
\textbf{Proof.}
This is based on standard techniques in the calculus
of variations and we quote \cite{MT2019} where relevant arguments can
be found.\quad $\Box$

\begin{proposition}\label{prop:3.3}
Under the above assumptions, we have:
$$
\lim_{\lambda\rightarrow 0}\frac{T_\lambda-T}{\lambda}
=-\frac{w_2(T)}{z_2^\prime (T)}
$$
if $z_2^\prime (T) \neq 0$.
\end{proposition}

\noindent
\textbf{Proof.}
Clearly $\nabla(g+\lambda r)\neq 0$ on $G_\lambda$ if
$|\lambda|$ small. Then, by the perturbed variant of
(\ref{2.9})-(\ref{2.11}) it yields
$|z_1^{\lambda \prime}(T_\lambda)|
+ |z_2^{\lambda \prime}(T_\lambda)| >0$ and,
similarly $|z_1^\prime (T)| + |z_2^\prime (T))| >0$, due to
(\ref{2.8}).
We choose here $z_2^\prime (T)\neq 0$ and, consequently,
$z_2^{\lambda \prime}(T_\lambda)\neq 0$, for $\lambda$ ``small''.
Then ${z}_2^\lambda$ is invertible on some interval
$[T-\alpha,T+\beta]$ with $\alpha, \beta >0$, small, not depending on $\lambda$, 
(and similarly around 0 due to the periodicity property).

This is due to $\mathbf{z}_{\lambda}\rightarrow \mathbf{z}$
in $\mathcal{C}^1([0,2T])^2$ and $T_\lambda \rightarrow T$.
We have $\mathbf{z}_{\lambda}(T_\lambda)=\mathbf{x}_0$ and
it yields:
\begin{equation}\label{3.11}
T_\lambda=(z_2^\lambda)^{-1}(x_0^2).
\end{equation}
We denote $x_0^\lambda=z_2(T_\lambda)\rightarrow x_0^2$
as $\lambda\rightarrow 0$. We may write
\begin{equation}\label{3.12}
\frac{T_\lambda-T}{\lambda}=
\frac{(z_2^\lambda)^{-1}(x_0^2)-(z_2)^{-1}(x_0^2)}{\lambda}
=\frac{(z_2)^{-1}(x_0^\lambda)-(z_2)^{-1}(x_0^2)}{\lambda}.
\end{equation}

By (\ref{3.11}), (\ref{3.12}) we get
$$
\frac{T_\lambda-T}{\lambda}=
\frac{(z_2)^{-1}(x_0^\lambda)-(z_2)^{-1}(x_0^2)}{x_0^\lambda-x_0^2}
\frac{z_2(T_\lambda)-z_2^\lambda(T_\lambda)}{\lambda}.
$$
Passing to the limit in the above relation and using
Proposition \ref{prop:3.2}, we end the proof.
\quad $\Box$

\begin{remark}\label{rem:3.2}
If $z_1^\prime (T)\neq 0$, the limit is
$-\frac{w_1(T)}{z_1^\prime (T)}$. In general, we denote by
$\theta(g,r)$ this limit. The last condition in Proposition \ref{prop:3.3} is a
consequence of (\ref{2.8}).
\end{remark}

To study the differentiability properties of the penalized
cost function (\ref{3.1}), we also assume
$f\in W^{1,p}(D)$, $\partial D$ is in $\mathcal{C}^{2,1}$
and $\mathcal{F}\subset \mathcal{C}^2(\overline{D})$.
We get that $g_+^2\in W^{1,\infty}(D)$ and
$g_+^2u\in W^{1,p}(D)$ if $u\in W^{1,p}(D)$ and the solution of
(\ref{3.2}), (\ref{3.3}) satisfies $y\in  W^{3,p}(D)
\subset \mathcal{C}^2(\overline{D})$.

\begin{proposition}\label{prop:3.4}
Under the above conditions, assume that 
$J(\mathbf{x},\cdot)$ is in $\mathcal{C}^1(\mathbb{R})$ and $j(\cdot,\cdot)$ is in
$\mathcal{C}^1(\mathbb{R}^3)$. Then, the directional derivative
of (\ref{3.1}), in the direction
$[v,r]\in W^{1,p}(D) \times \mathcal{F}$, is given by:
\begin{eqnarray}
&&
\theta(g,r)\left[
j(\mathbf{x}_0,y(\mathbf{x}_0))
+\left|\frac{\partial y}{\partial \mathbf{n}}(\mathbf{x}_0) \right|^2
\right] | \nabla g(\mathbf{x}_0)|
+\int_E \partial_2 J(\mathbf{x},y(\mathbf{x})) q(\mathbf{x})d\mathbf{x}
\nonumber\\
&+&
\int_0^T
\nabla_1 j\left(\mathbf{z}(t),y(\mathbf{z}(t))\right)
\cdot \mathbf{w}(t)|\mathbf{z}^\prime(t)| dt
\nonumber\\
&+&
\int_0^T
\partial_2 j\left(\mathbf{z}(t),y(\mathbf{z}(t))\right)
\left[
\nabla y(\mathbf{z}(t))\cdot \mathbf{w}(t)
+q(\mathbf{z}(t))
\right]
|\mathbf{z}^\prime(t)| dt
\nonumber\\
&+&
\int_0^T
j\left(\mathbf{z}(t),y(\mathbf{z}(t))\right)
\frac{\mathbf{z}^\prime(t)\cdot \mathbf{w}^\prime(t)}
{|\mathbf{z}^\prime(t)|}
dt
\nonumber\\
&+&\frac{2}{\epsilon}
\int_0^T
\nabla y(\mathbf{z}(t))\cdot
\frac{\nabla g(\mathbf{z}(t))}{|\nabla g(\mathbf{z}(t))|^2}
\nabla r(\mathbf{z}(t))\cdot\nabla y(\mathbf{z}(t))
|\mathbf{z}^\prime(t)|
dt
\nonumber\\
&+&\frac{2}{\epsilon}
\int_0^T
\nabla y(\mathbf{z}(t))\cdot
\frac{\nabla g(\mathbf{z}(t))}{|\nabla g(\mathbf{z}(t))|}
\left[
\left(H\,y(\mathbf{z}(t))\right)\mathbf{w}(t)
+ \nabla q(\mathbf{z}(t))
\right]
\cdot
\frac{\nabla g(\mathbf{z}(t))}{|\nabla g(\mathbf{z}(t))|}
|\mathbf{z}^\prime(t)|
dt
\nonumber\\
&+&\frac{2}{\epsilon}
\int_0^T
\nabla y(\mathbf{z}(t))\cdot
\frac{\nabla g(\mathbf{z}(t))}{|\nabla g(\mathbf{z}(t))|}
\nabla y(\mathbf{z}(t))\cdot
\left[
\frac{\left(H\,g(\mathbf{z}(t))\right)\mathbf{w}(t)}
{|\nabla g(\mathbf{z}(t))|}
\right.  
\nonumber\\
&&-\left.
\frac{\nabla g(\mathbf{z}(t))}{|\nabla g(\mathbf{z}(t))|^3}
\left(
\nabla g(\mathbf{z}(t))\cdot \nabla r(\mathbf{z}(t))
+\nabla g(\mathbf{z}(t))
\left(H\,g(\mathbf{z}(t))\right)\mathbf{w}(t)
\right)
\right]|\mathbf{z}^\prime(t)|
dt
\nonumber\\
&+&\frac{1}{\epsilon}
\int_0^T
\left[
\nabla y(\mathbf{z}(t))\cdot
\frac{\nabla g(\mathbf{z}(t))}{|\nabla g(\mathbf{z}(t))|}  
\right]^2
\frac{\mathbf{z}^\prime(t)\cdot \mathbf{w}^\prime(t)}
{|\mathbf{z}^\prime(t)|}
dt .
\label{3.13}
\end{eqnarray}

\end{proposition}

The notations are explained in the proof.

\noindent
\textbf{Proof.}
We compute
\begin{eqnarray*}
&&
\lim_{\lambda\rightarrow 0} \frac{1}{\lambda}
\left\{
\int_E J(\mathbf{x},y_\lambda(\mathbf{x})) d\mathbf{x}
+\int_0^{T_\lambda}
j\left(\mathbf{z}_\lambda(t),y_\lambda(\mathbf{z}_\lambda(t))\right)
|\mathbf{z}_\lambda^\prime(t)|
dt
\right.
\nonumber\\
&&
+\frac{1}{\epsilon}
\int_0^{T_\lambda}
\left[
\nabla y_\lambda(\mathbf{z}_\lambda(t))\cdot
\frac{\nabla (g+\lambda r)(\mathbf{z}_\lambda(t))}
{|\nabla (g+\lambda r)(\mathbf{z}_\lambda(t))|}
\right]^2     
|\mathbf{z}_\lambda^\prime(t)|
dt
-\int_E J(\mathbf{x},y(\mathbf{x})) d\mathbf{x}
\nonumber\\
&&
\left.
-\int_0^T
j\left(\mathbf{z}(t),y(\mathbf{z}(t))\right)
|\mathbf{z}^\prime(t)|
dt
-\frac{1}{\epsilon}
\int_0^T
\left[
\nabla y(\mathbf{z}(t))\cdot
\frac{\nabla g(\mathbf{z}(t))}{|\nabla g(\mathbf{z}(t))|}  
\right]^2  
|\mathbf{z}^\prime(t)|
dt
\right\} .
\end{eqnarray*}

Applying Proposition \ref{prop:3.2}, (\ref{3.6}),
(\ref{3.7}), and the differentiability hypotheses on $J$,
$j$, we get:
\begin{equation}\label{3.14}
\frac{1}{\lambda}
\left[
\int_E J(\mathbf{x},y_\lambda(\mathbf{x})) d\mathbf{x}
-\int_E J(\mathbf{x},y(\mathbf{x})) d\mathbf{x}
\right]
\rightarrow
\int_E
\partial_2 J(\mathbf{x},y(\mathbf{x}))q(\mathbf{x})
d\mathbf{x}.
\end{equation}

We discuss now the term:
\begin{eqnarray}
&&\frac{1}{\lambda}
\int_T^{T_\lambda}
j\left(\mathbf{z}_\lambda(t),y_\lambda(\mathbf{z}_\lambda(t))\right)
|\mathbf{z}_\lambda^\prime(t)|
dt
=\frac{T_\lambda-T}{\lambda}
j\left(\mathbf{z}_\lambda(\tau_\lambda),
y_\lambda(\mathbf{z}_\lambda(\tau_\lambda))\right)
|\mathbf{z}_\lambda^\prime(\tau_\lambda)|
\nonumber\\
&&\rightarrow
\theta(g,r)j(\mathbf{x}_0,y(\mathbf{x}_0))
|\mathbf{z}^\prime(T)|
=\theta(g,r)j(\mathbf{x}_0,y(\mathbf{x}_0))
|\nabla g(\mathbf{x}_0) |,
\label{3.15}
\end{eqnarray}
due to (\ref{2.9})-(\ref{2.11}) and Remark \ref{rem:3.2}.
Here $\tau_\lambda$ is some intermediary point in the
interval $[T, T_\lambda]$, depending on $\lambda$, $g$, $r$,
$j$, etc. We also use Thm. \ref{theo:2.1} and
$T_\lambda \rightarrow T$.

Similarly, we consider the term:
\begin{eqnarray}
&&\frac{1}{\lambda}
\int_T^{T_\lambda}
\left[
\nabla y_\lambda(\mathbf{z}_\lambda(t))\cdot
\frac{\nabla (g+\lambda r)(\mathbf{z}_\lambda(t))}
{|\nabla (g+\lambda r)(\mathbf{z}_\lambda(t))|}
\right]^2     
|\mathbf{z}_\lambda^\prime(t)|
dt
\nonumber\\
&\rightarrow&
\theta(g,r)\left[
  \nabla y(\mathbf{x}_0)\cdot
\frac{\nabla g(\mathbf{x}_0)}
{|\nabla g(\mathbf{x}_0)|}  
  \right]^2
|\nabla g(\mathbf{x}_0) |
=\theta(g,r)
\left| \frac{\partial y}
       {\partial \mathbf{n}}(\mathbf{x}_0)\right|^2
|\nabla g(\mathbf{x}_0) | .
\label{3.16}
\end{eqnarray}
In the last two limits, the regularity properties of
$y$, $\mathbf{z}$, $y_\lambda$, $\mathbf{z}_\lambda$ also play
a key role.

Next, we investigate the last term:
\begin{eqnarray*}
&&
\frac{1}{\lambda}
\left\{
\int_0^{T}
j\left(\mathbf{z}_\lambda(t),y_\lambda(\mathbf{z}_\lambda(t))\right)
|\mathbf{z}_\lambda^\prime(t)|
dt
\right.
\nonumber\\
&&
+
\frac{1}{\epsilon}
\int_0^{T}
\left[
\nabla y_\lambda(\mathbf{z}_\lambda(t))\cdot
\frac{\nabla (g+\lambda r)(\mathbf{z}_\lambda(t))}
{|\nabla (g+\lambda r)(\mathbf{z}_\lambda(t))|}
\right]^2 
|\mathbf{z}_\lambda^\prime(t)|
dt
\nonumber\\
&&
\left.
-\int_0^{T}
j\left(\mathbf{z}(t),y(\mathbf{z}(t))\right)
|\mathbf{z}^\prime(t)|
dt
-\frac{1}{\epsilon}
\int_0^{T}
\left[
\nabla y(\mathbf{z}(t))\cdot
\frac{\nabla g(\mathbf{z}(t))}
{|\nabla g(\mathbf{z}(t))|}
\right]^2
|\mathbf{z}^\prime(t)|
dt
\right\}
\end{eqnarray*}
Clearly, the terms containing $j(\cdot,\cdot)$ give
the limit:
\begin{eqnarray}
&&\int_0^{T}
\left[ \nabla_1 j\left(\mathbf{z}(t),y(\mathbf{z}(t))\right)
\cdot \mathbf{w}(t) 
+\partial_2 j\left(\mathbf{z}(t),y(\mathbf{z}(t))\right)
\nabla y(\mathbf{z}(t)) \cdot \mathbf{w}(t)
\right]
|\mathbf{z}^\prime(t)|
dt
\nonumber\\
&+&
\int_0^{T}
\left[
\partial_2 j\left(\mathbf{z}(t),y(\mathbf{z}(t))\right)
q(\mathbf{z}(t))
|\mathbf{z}^\prime(t)| 
+j\left(\mathbf{z}(t),y(\mathbf{z}(t))\right)
\frac{\mathbf{z}^\prime(t)\cdot \mathbf{w}^\prime(t)}
     {|\mathbf{z}^\prime(t)|}
\right]  
dt
\label{3.17}  
\end{eqnarray}
where $\nabla_1 j$ is the gradient of $j(\cdot,\cdot)$
with respect to the two components of $\mathbf{z}$, and
$\partial_2 j$ is the partial derivative with respect to $y$,
other quantities are defined in (\ref{3.6})-(\ref{3.10}).

Let us consider now the two terms corresponding to the
penalization of Neumann boundary condition.
We intercalate advantageous terms and we compute step by step:
\begin{eqnarray}
&&\frac{1}{\lambda}\int_0^{T}
\left\{
\left[ 
\nabla y_\lambda(\mathbf{z}_\lambda(t))\cdot
\frac{\nabla (g+\lambda r)(\mathbf{z}_\lambda(t))}
{|\nabla (g+\lambda r)(\mathbf{z}_\lambda(t))|}
\right]^2
-\left[
\nabla y(\mathbf{z}(t))\cdot
\frac{\nabla g(\mathbf{z}(t))}
{|\nabla g(\mathbf{z}(t))|}
\right]^2
\right\}
|\mathbf{z}_\lambda^\prime(t)| dt
\nonumber\\
&=&\frac{1}{\lambda}
\int_0^{T}
S\left[
\nabla y_\lambda(\mathbf{z}_\lambda(t))\cdot
\frac{\nabla (g+\lambda r)(\mathbf{z}_\lambda(t))}
 {|\nabla (g+\lambda r)(\mathbf{z}_\lambda(t))|}
-\nabla y(\mathbf{z}(t))\cdot
\frac{\nabla g(\mathbf{z}(t))}
{|\nabla g(\mathbf{z}(t))|}     
\right]  
|\mathbf{z}_\lambda^\prime(t)| dt
\nonumber\\
&=&
\int_0^{T}
S\frac{\nabla r(\mathbf{z}_\lambda(t))}
{|\nabla (g+\lambda r)(\mathbf{z}_\lambda(t))|}
\cdot
\nabla y_\lambda(\mathbf{z}_\lambda(t))
|\mathbf{z}_\lambda^\prime(t)| dt
\nonumber\\
&&
+\int_0^{T}
S\frac{\nabla y_\lambda(\mathbf{z}_\lambda(t))
  -\nabla y(\mathbf{z}(t))}{\lambda}
\cdot
\frac{\nabla g(\mathbf{z}(t))}
{|\nabla g(\mathbf{z}(t))|}
|\mathbf{z}_\lambda^\prime(t)| dt
\nonumber\\
&&
+\frac{1}{\lambda}
\int_0^{T}
S\left[\nabla
y_\lambda(\mathbf{z}_\lambda(t))\cdot
\frac{\nabla g(\mathbf{z}_\lambda(t))}
 {|\nabla (g+\lambda r)(\mathbf{z}_\lambda(t))|}
 -\nabla y_\lambda(\mathbf{z_\lambda}(t))\cdot
\frac{\nabla g(\mathbf{z}(t))}
 {|\nabla (g)(\mathbf{z}(t))|}
 \right]
|\mathbf{z}_\lambda^\prime(t)| dt
\nonumber\\
&=&I+II+III
\label{3.18}  
\end{eqnarray}
where $S$ is the sum
$$
\nabla y_\lambda(\mathbf{z}_\lambda(t))\cdot
\frac{\nabla (g+\lambda r)(\mathbf{z}_\lambda(t))}
 {|\nabla (g+\lambda r)(\mathbf{z}_\lambda(t))|}
+\nabla y(\mathbf{z}(t))\cdot
\frac{\nabla g(\mathbf{z}(t))}
{|\nabla g(\mathbf{z}(t))|} .
$$

We have:
\begin{eqnarray*}
\lim_{\lambda \rightarrow 0} I & = &
2\int_0^T
\nabla y(\mathbf{z}(t))\cdot
\frac{\nabla g(\mathbf{z}(t))}{|\nabla g(\mathbf{z}(t))|}
\frac{\nabla r(\mathbf{z}(t))}{|\nabla g(\mathbf{z}(t))|}
\cdot\nabla y(\mathbf{z}(t))
|\mathbf{z}^\prime(t)|
dt
\\
\lim_{\lambda \rightarrow 0} II & = &
2\int_0^T
\nabla y(\mathbf{z}(t))\cdot
\frac{\nabla g(\mathbf{z}(t))}{|\nabla g(\mathbf{z}(t))|}
\left[
H\,y(\mathbf{z}(t)) + \nabla q(\mathbf{z}(t))
\right]
\cdot
\frac{\nabla g(\mathbf{z}(t))}{|\nabla g(\mathbf{z}(t))|}
|\mathbf{z}^\prime(t)|
dt,
\end{eqnarray*}
where $H\,y$ is the Hessian matrix of
$y\in\mathcal{C}^2(\overline{D})$.

Concerning part $III$, we get:
\begin{eqnarray*}
\lim_{\lambda \rightarrow 0} III & = &
2\int_0^T
\nabla y(\mathbf{z}(t))\cdot
\frac{\nabla g(\mathbf{z}(t))}{|\nabla g(\mathbf{z}(t))|}
|\mathbf{z}^\prime(t)|
\nabla y(\mathbf{z}(t))\cdot
\left[
\frac{\left(H\,g(\mathbf{z}(t))\right)\mathbf{w}(t)}
 {|\nabla g(\mathbf{z}(t))|}
 \right.
\\
&&
\left.
-\frac{\nabla g(\mathbf{z}(t))}{|\nabla g(\mathbf{z}(t))|^2}
\left(
\frac{\nabla g(\mathbf{z}(t))\cdot \nabla r(\mathbf{z}(t))}
     {|\nabla g(\mathbf{z}(t))|}
     +\frac{\nabla g(\mathbf{z}(t))\cdot
     \left(H\,g(\mathbf{z}(t))\right)\mathbf{w}(t)}
     {|\nabla g(\mathbf{z}(t))|}
\right)
\right]
dt
\\
&=&
2\int_0^T
\nabla y(\mathbf{z}(t))\cdot
\frac{\nabla g(\mathbf{z}(t))}{|\nabla g(\mathbf{z}(t))|}
\nabla y(\mathbf{z}(t))\cdot
\left[
\frac{\left(H\,g(\mathbf{z}(t))\right)\mathbf{w}(t)}
{|\nabla g(\mathbf{z}(t))|}
\right.  
\\
&&
-\left.
\frac{\nabla g(\mathbf{z}(t))}{|\nabla g(\mathbf{z}(t))|^3}
\left(
\nabla g(\mathbf{z}(t))\cdot \nabla r(\mathbf{z}(t))
+\nabla g(\mathbf{z}(t))
\left(H\,g(\mathbf{z}(t))\right)\mathbf{w}(t)
\right)
\right]|\mathbf{z}^\prime(t)|
dt
\end{eqnarray*}

Finally, the term
\begin{eqnarray}
&&\int_0^T
\left[
\nabla y(\mathbf{z}(t))\cdot
\frac{\nabla g(\mathbf{z}(t))}{|\nabla g(\mathbf{z}(t))|}  
\right]^2
\frac{|\mathbf{z}_\lambda^\prime(t)| - |\mathbf{z}^\prime(t)|}
{\lambda}
dt
\nonumber\\  
&\rightarrow &  
\int_0^T
\left[
\nabla y(\mathbf{z}(t))\cdot
\frac{\nabla g(\mathbf{z}(t))}{|\nabla g(\mathbf{z}(t))|}  
\right]^2
\frac{\mathbf{z}^\prime(t)\cdot \mathbf{w}^\prime(t)}
{|\mathbf{z}^\prime(t)|}
dt
\label{3.19}  
\end{eqnarray}

Summing up relations (\ref{3.14})-(\ref{3.19}), we
finish the proof of (\ref{3.13}).
\quad$\Box$

\section{Finite element descent directions}
\setcounter{equation}{0}

We use the piecewise cubic finite element $\mathbb{P}_3$ in
$\mathcal{T}_h$ a triangulation of $D$.
We define 
$$
\mathbb{W}_h=\{ \varphi_h\in \mathcal{C}(\overline{D});
\ {\varphi_h}_{|T} \in \mathbb{P}_3(T),\ \forall T \in \mathcal{T}_h  \}
$$
of dimension $n=card(I)$ ($I$ the set of nodes in
$\mathcal{T}_h$) and
$$
\mathbb{V}_h=\{ \varphi_h\in \mathbb{W}_h;
\ \varphi_h = 0\hbox{ on }\partial D \},
$$
of dimension $n_0=card(I_0)$ ($I_0$ the set of nodes in
$\mathcal{T}_h$, outside $\partial D$) which are finite element approximations
of Hilbert spaces $\mathbb{W}=H^1(D)$, $\mathbb{V}=H^1_0(D)$, respectively.

The parametrization function $g$ is approached by the finite element function
$g_h\in \mathbb{W}_h$, $g_h(\mathbf{x})=\sum_{i\in I} G_i \phi_i(\mathbf{x})$
where $G=(G_i)_{i\in I} \in \mathbb{R}^n$ is a real vector and $\phi_i$ is
the basis in $\mathbb{W}_h$. Similarly, we denote
$u_h\in \mathbb{W}_h$, $y_h\in \mathbb{V}_h$ and the associated vectors
$U=(U_i)_{i\in I}\in\mathbb{R}^n$ and $Y=(Y_j)_{j\in I_0}\in\mathbb{R}^{n_0}$ for the
discretization of the control, respectively the state.
For the control term $u_h$, one can also employ lower order finite elements,
like continuous piecewise linear $\mathbb{P}_1$
or piecewise constant $\mathbb{P}_0$. See \cite{Ciarlet2002}, \cite{Raviart2004}
for a discussion of finite element spaces.

Here, we  consider (\ref{2.1}) with non homogeneous
boundary condition $\frac{\partial y_\Omega}{\partial \mathbf{n}}=\delta$
on $\partial \Omega$, with $\delta$ some given function in $H^1 (D)$.
The objective function (\ref{3.1}) is taken of the form
\begin{eqnarray}
\min_{g,u} \mathcal{J}(g,u)&=&
\left\{
\int_{E} J\left(\mathbf{x},y(\mathbf{x})\right)
d\mathbf{x}
+
\int_{I_g}
j\left(\mathbf{z}(t), y(\mathbf{z}(t))\right)
|\mathbf{z}^\prime(t)| dt
\right.
\nonumber\\
&&
\left.
+\frac{1}{\epsilon}
\int_{I_g}
\left[
  \nabla y(\mathbf{z}(t))\cdot
  \frac{\nabla g(\mathbf{z}(t)) }{
    |\nabla g(\mathbf{z}(t)) |}
  -\delta(\mathbf{z}(t))
\right]^2
|\mathbf{z}^\prime(t)| dt
\right\} .
\label{5.1}
\end{eqnarray}

We denote the first term of (\ref{5.1}) by
$$
t_1=\int_{E} J\left(\mathbf{x},y(\mathbf{x})\right)
d\mathbf{x}.
$$
The second and the third terms of (\ref{5.1}) can be rewritten as integrals on
$\partial \Omega_g$, more precisely
\begin{eqnarray*}
t_2 &=& \int_{\partial \Omega_g}
j\left(s, y(s)\right) ds \\
t_3 &=& \frac{1}{\epsilon}
\int_{\partial \Omega_g}
\left[
  \nabla y(s)\cdot
  \frac{\nabla g(s) }{
    |\nabla g(s) |}
  -\delta(s)
  \right]^2 ds .
\end{eqnarray*} 
We employ the software FreeFem++, \cite{freefem++} and these terms
can be computed with the command 
\texttt{int1d(Th,levelset=gh)(\dots)}.

We use the general descent direction method
$$
(G^{k+1},U^{k+1})=(G^k,U^k)+\lambda_k (R^k,V^k),
$$
where $\lambda_k >0$ is obtained via some line search
$$
\lambda_k \in \arg\min_{\lambda >0}
\mathcal{J}\left((G^k,U^k)+\lambda (R^k,V^k)\right)
$$
and $(R^k,V^k)$ is a descent direction, i.e. $d\mathcal{J}_{(G^k,U^k)}(R^k,V^k) <0$.
For $E \neq \emptyset$, a projection is necessary in order to get (\ref{2.4}).
The algorithm stops if
$| \mathcal{J}(G^{k+1},U^{k+1}) - \mathcal{J}(G^k,U^k)| < tol$ or
$d\mathcal{J}_{(G^k,U^k)}(R^k,V^k)=0$. Other choices are possible, see \cite{Ciarlet2018}
for details on such algorithms.

Since the approximating state system (\ref{3.2}), (\ref{3.3}) is similar
to \cite{MT2019}, we apply here a similar discretization technique of
the gradient (\ref{3.13}). In the following, we shall
use descent directions based on the discrete simplified adjoint system:
find $p_h\in \mathbb{V}_h$ such that
\begin{eqnarray}
&&\int_D \nabla \varphi_h \cdot \nabla p_h d\mathbf{x}
+\int_D \varphi_h p_h d\mathbf{x}
= 
\int_{E} \partial_2 J\left(\mathbf{x},y_h(\mathbf{x})\right)
\varphi_h(\mathbf{x})
d\mathbf{x}
\nonumber\\
&&
+\int_{\partial \Omega_{g_h}} \partial_2 j\left(s, y_h(s)\right)
\varphi_h(s) ds
\nonumber\\
&&
+\frac{2}{\epsilon}
\int_{\partial \Omega_{g_h}}
\left(
  \nabla y_h(s)\cdot
  \frac{\nabla g_h(s) }{
    |\nabla g_h(s) |}
  -\delta_h(s)
  \right)
\nabla \varphi_h(s) \cdot
  \frac{\nabla g_h(s) }{
    |\nabla g_h(s) |}
  ds
\label{5.2}
\end{eqnarray}
for all $\varphi_h \in \mathbb{V}_h$.
In the right hand side of (\ref{5.2}) appear just the terms multiplying $q$
in the gradient (\ref{3.13}) and $\delta_h(s)$ is a continuous piecewise
linear $\mathbb{P}_1$ discretization of $\delta(s)$ in $D$.

\begin{proposition}\label{prop:5.1}
Given $g_h,u_h\in\mathbb{W}_h$ and the variations $r_h,v_h\in\mathbb{W}_h$,
let $y_h\in\mathbb{V}_h$ be the
finite element solution of
(\ref{3.2}), (\ref{3.3}), let
$q_h\in\mathbb{V}_h$ be the finite element solution of (\ref{3.6}),
(\ref{3.7}) depending in $r_h,\ v_h$
and let $p_h\in\mathbb{V}_h$ be the solution of (\ref{5.2}).
Then
\begin{eqnarray}
&&\int_{E} \partial_2 J\left(\mathbf{x},y_h(\mathbf{x})\right)
q_h(\mathbf{x})
d\mathbf{x}
+\int_{\partial \Omega_{g_h}} \partial_2 j\left(s, y_h(s)\right)
q_h(s) ds
\nonumber\\
&&
+\frac{2}{\epsilon}
\int_{\partial \Omega_{g_h}}
\left(
  \nabla y_h(s)\cdot
  \frac{\nabla g_h(s) }{
    |\nabla g_h(s) |}
  -\delta_h(s)
  \right)
\nabla q_h(s) \cdot
  \frac{\nabla g_h(s) }{
    |\nabla g_h(s) |}
  ds
\leq 0  
\label{5.3}  
\end{eqnarray}
if we choose:\\
i) $r_h=-p_hu_h$ and $v_h=-p_h$ or\\
ii) $r_h=-\widetilde{d}_h$ and $v_h=-p_h$ where
$\widetilde{d}_h \in\mathbb{W}_h$ is the solution of
\begin{eqnarray}
&&\int_D \nabla \widetilde{d}_h \cdot \nabla \varphi_h d\mathbf{x}
+\int_D \widetilde{d}_h \varphi_h d\mathbf{x}
= 
\int_D  2(g_h)_+ u_h p_h
\varphi_h
d\mathbf{x}
\label{5.4}
\end{eqnarray}
for all $\varphi_h \in \mathbb{W}_h$.
\end{proposition}  

\noindent
\textbf{Proof.}
Putting $\varphi_h=q_h$ in (\ref{5.2}) and multiplying (\ref{3.6}) by $p_h$,
integrating by parts over $D$ and using (\ref{3.7}),
we get that the left hand side of (\ref{5.3}) is equal to:
\begin{eqnarray*}
  \int_D (g_h)_+^2 v_h p_h d\mathbf{x}
  + \int_D 2(g_h)_+ u_h r_h p_h  d\mathbf{x}.  
\end{eqnarray*}
For $v_h=-p_h$, we have
$$
\int_D (g_h)_+^2 v_hp_hd\mathbf{x}
=-\int_D (g_h)_+^2 p_h^2 d\mathbf{x}
\leq 0.
$$
If $(g_h)_+ p_h$ is not null, then the above inequality is strict.\\
Case i). For $r_h=-p_h u_h$, we have
$$
\int_D  2(g_h)_+ u_h r_h p_h d\mathbf{x}
=-\int_D  2(g_h)_+ (u_h p_h)^2 d\mathbf{x}
\leq 0.
$$
Case ii). For $r_h=-\widetilde{d}_h$, we have
\begin{eqnarray*}
\int_D  2(g_h)_+ u_h r_h p_h d\mathbf{x}
&=&-\int_D  2(g_h)_+ u_h p_h \widetilde{d}_h d\mathbf{x}\\
&=&-\int_D \nabla \widetilde{d}_h \cdot \nabla \widetilde{d}_h d\mathbf{x}
-\int_D \widetilde{d}_h \widetilde{d}_h d\mathbf{x}
\leq 0.
\end{eqnarray*}
The second equality is obtained by putting
$\varphi_h =\widetilde{d}_h$ in (\ref{5.4}).
This ends the proof.
If $(g_h)_+ p_h$ is not null, then the inequality (\ref{5.3}) is strict.
\quad$\Box$

\begin{remark}\label{rem:5.1}
Due to the strong non convex character of the shape optimization problems,
the descent algorithms find just a local minimum point of the penalized problem,
in general. The penalization term may remain not null, that is the
constraint (\ref{2.2}) may be violated. However, the above methodology offers
a systematic and general approximation procedure that can be applied in many
examples and produces relevant results. Both topological and boundary
variations are performed simultaneously.
\end{remark}

\section{Numerical tests}
\setcounter{equation}{0}

\medskip

\textbf{Example 1.}

We choose
$D=]-3,3[\times ]-3,3[$,
$y_d(x_1,x_2)=x_1^2 + x_2^2 -1^2$, $f(\mathbf{x})=-4+y_d(\mathbf{x})$   
and the tracking type cost
$j(\mathbf{x})=\frac{1}{2}\left(y(\mathbf{x})-y_d(\mathbf{x})\right)^2$.
We fix $\delta=2$ for the non homogeneous Neumann
boundary condition.
We consider first the case $E=\emptyset$ and $J = 0$, with the numerical parameters:
$\epsilon=0.5$, the mesh of $D$ has 73786 triangles and 37254 vertices
and the tolerance parameter for the stopping test is $tol=10^{-6}$.

The initial domain is the disk of center $(0,0)$ and radius $2.5$
with a circular hole of center $(-1,-1)$ and radius $0.5$. The corresponding
$g_0(x_1,x_2)$ is given by
$$
\max\left(
 (x_1)^2 + (x_2)^2 -2.5^2,
-(x_1+1)^2 - (x_2+1)^2 +0.5^2
\right) .
$$
The initial guess for the control is $u_0=0$.

We use the descent direction
given by the Proposition \ref{prop:5.1}, case ii) and
the algorithm stops after 3 iterations.
For the stopping test, we have computed just the left hand side of (\ref{5.3}) and we
replaced  $d\mathcal{J}_{(G^k,U^k)}(R^k,V^k)=0$ by:
there are no smaller values than $\mathcal{J}(G^k,U^k)$ in the direction
$(R^k,V^k)$ for $\lambda  \in \{\rho^i; i\in \mathbb{N},\ 0 \leq i < 30 \}$,
with $\rho=0.8$.

\begin{figure}[ht]
\begin{center}
\includegraphics[width=4.5cm]{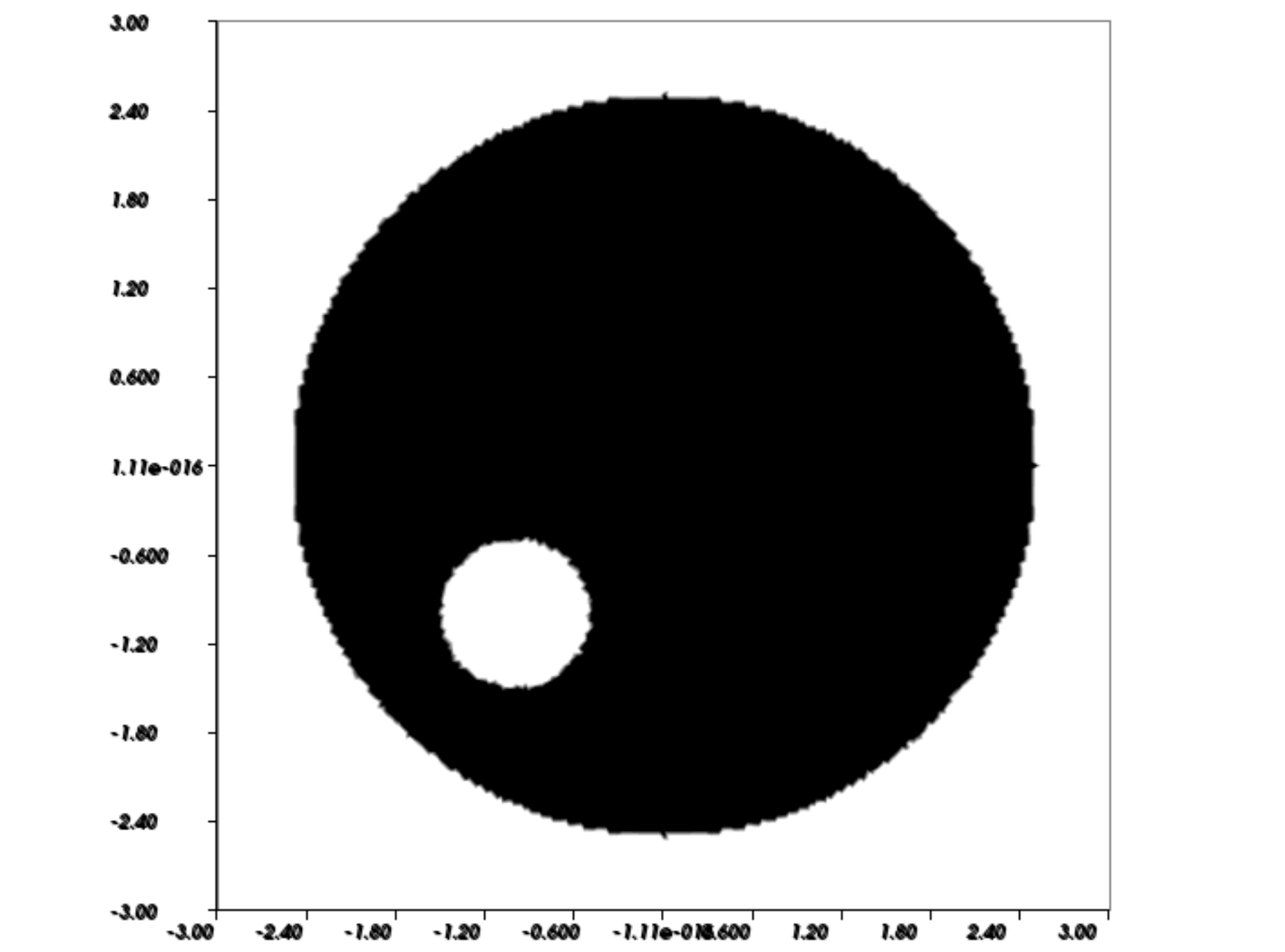}
\   
\includegraphics[width=4.5cm]{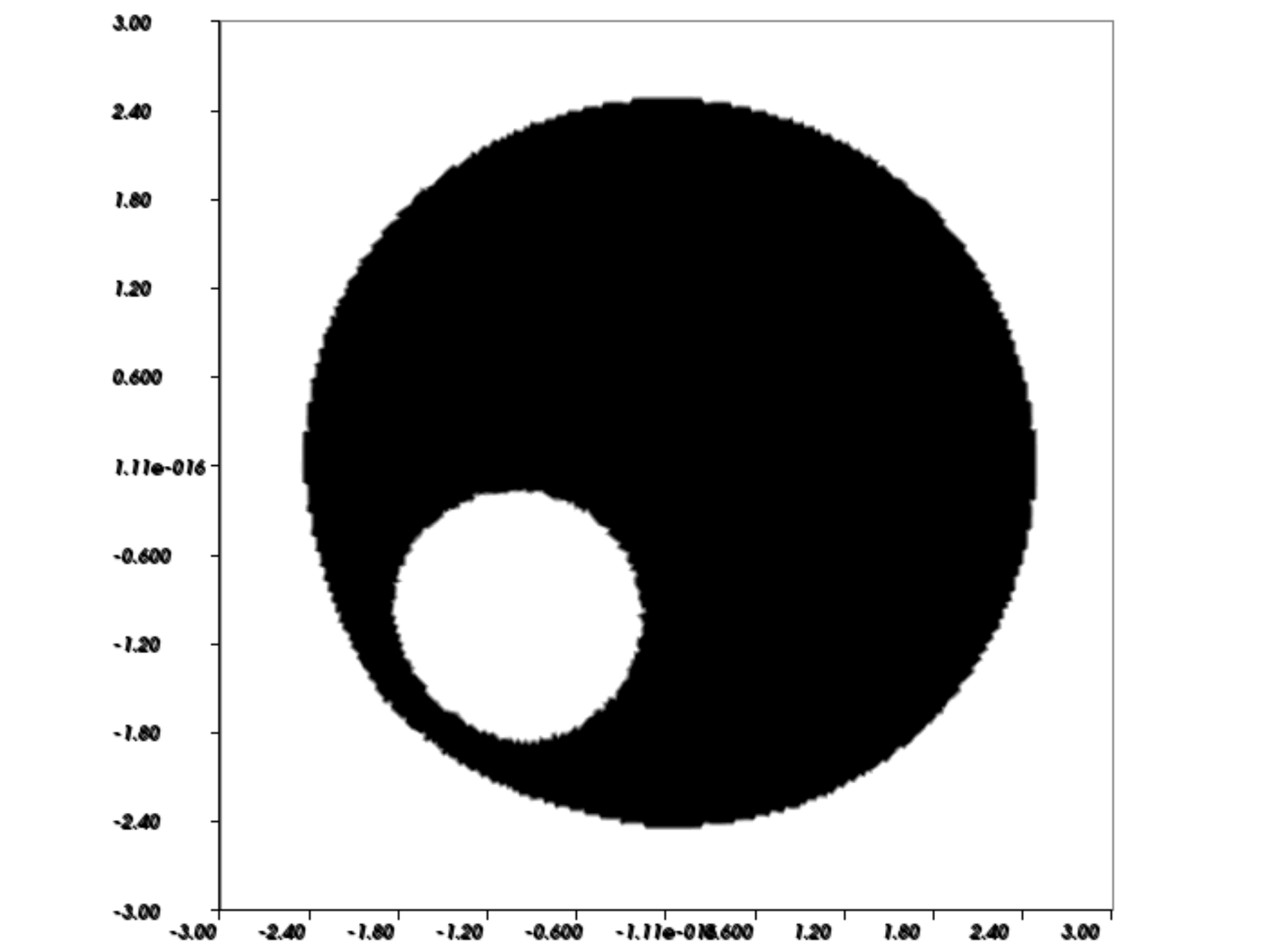}
\
\includegraphics[width=4.5cm]{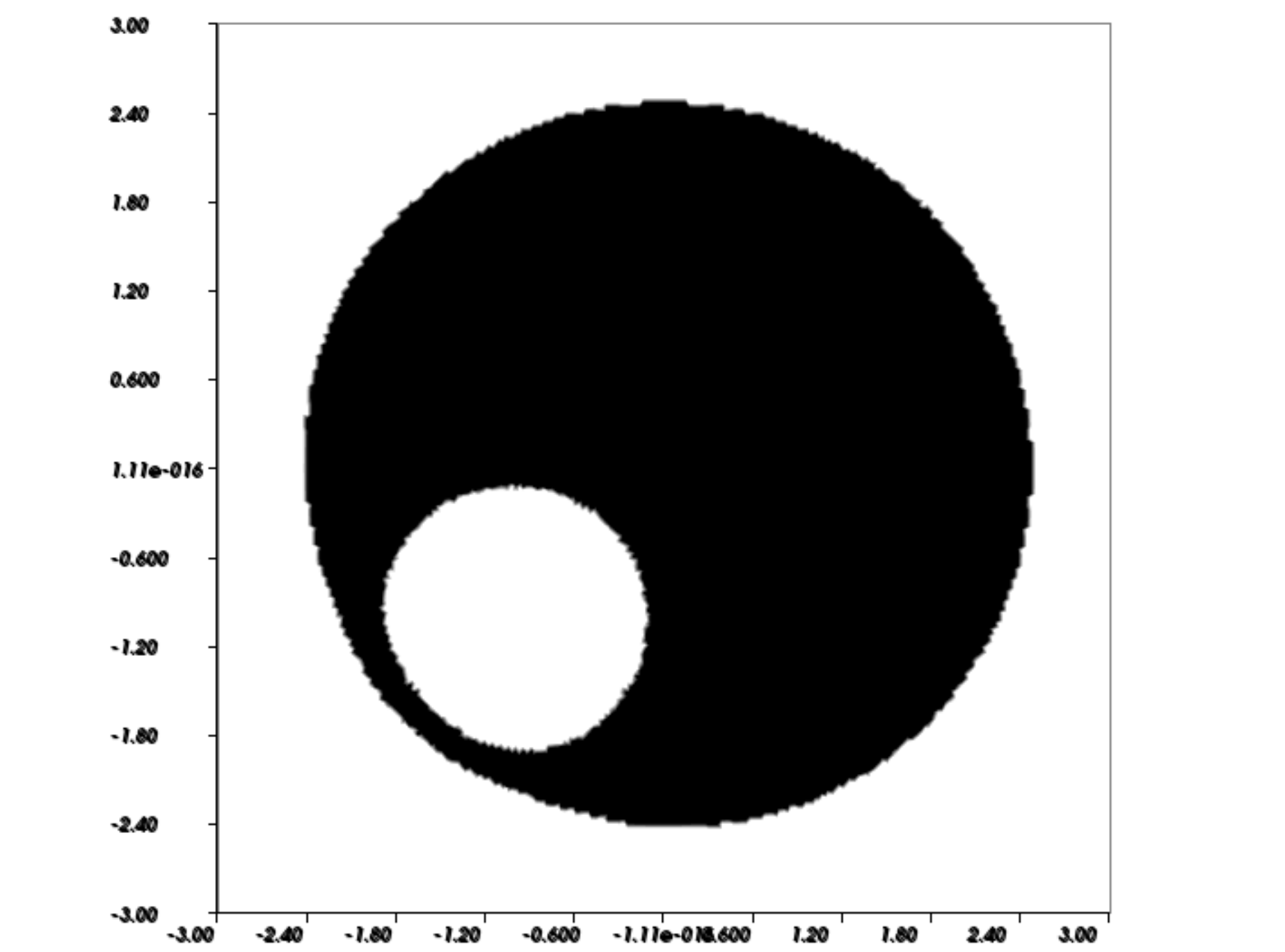}\\
\includegraphics[width=4.5cm]{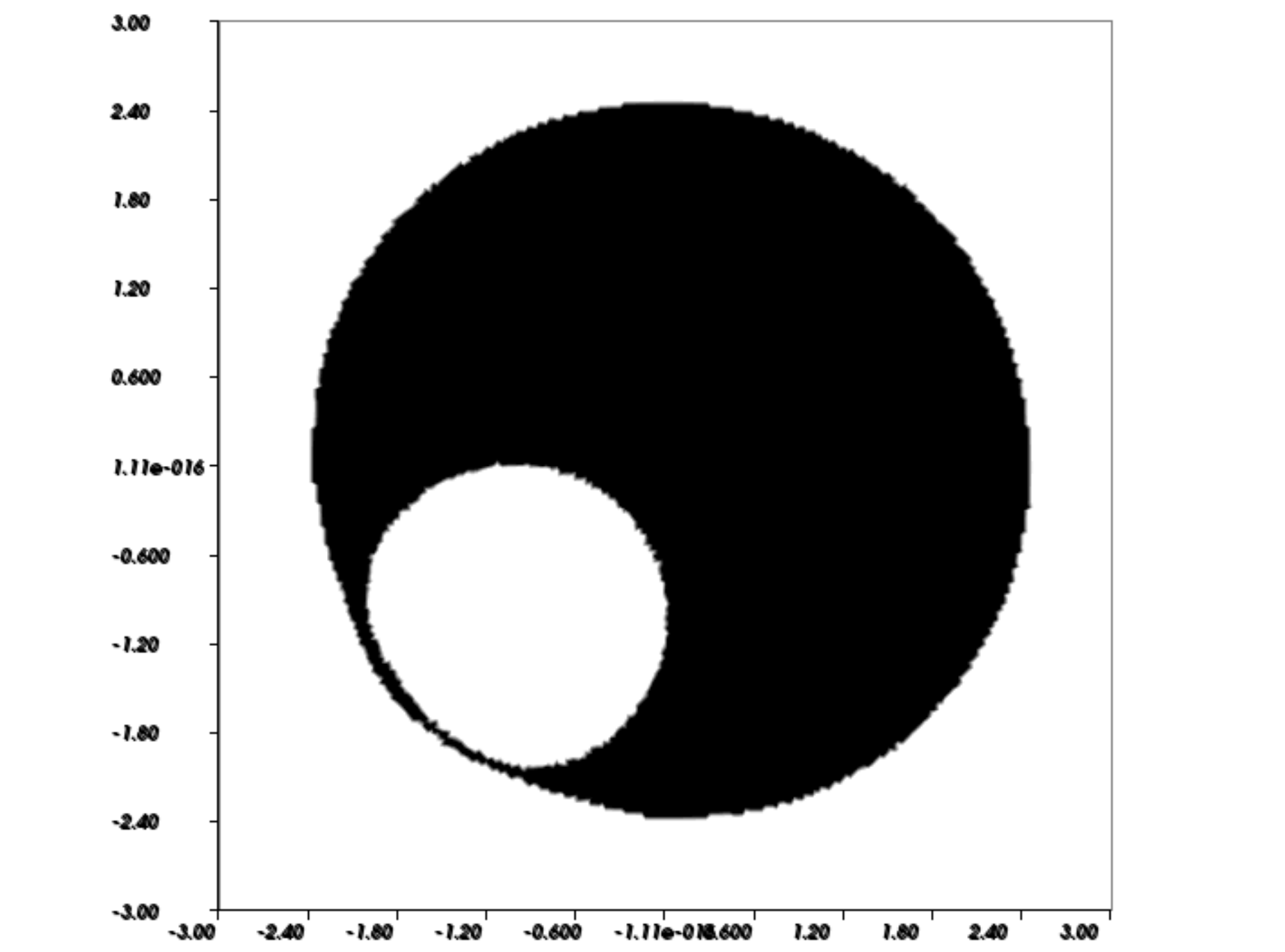}
\
\includegraphics[width=4.5cm]{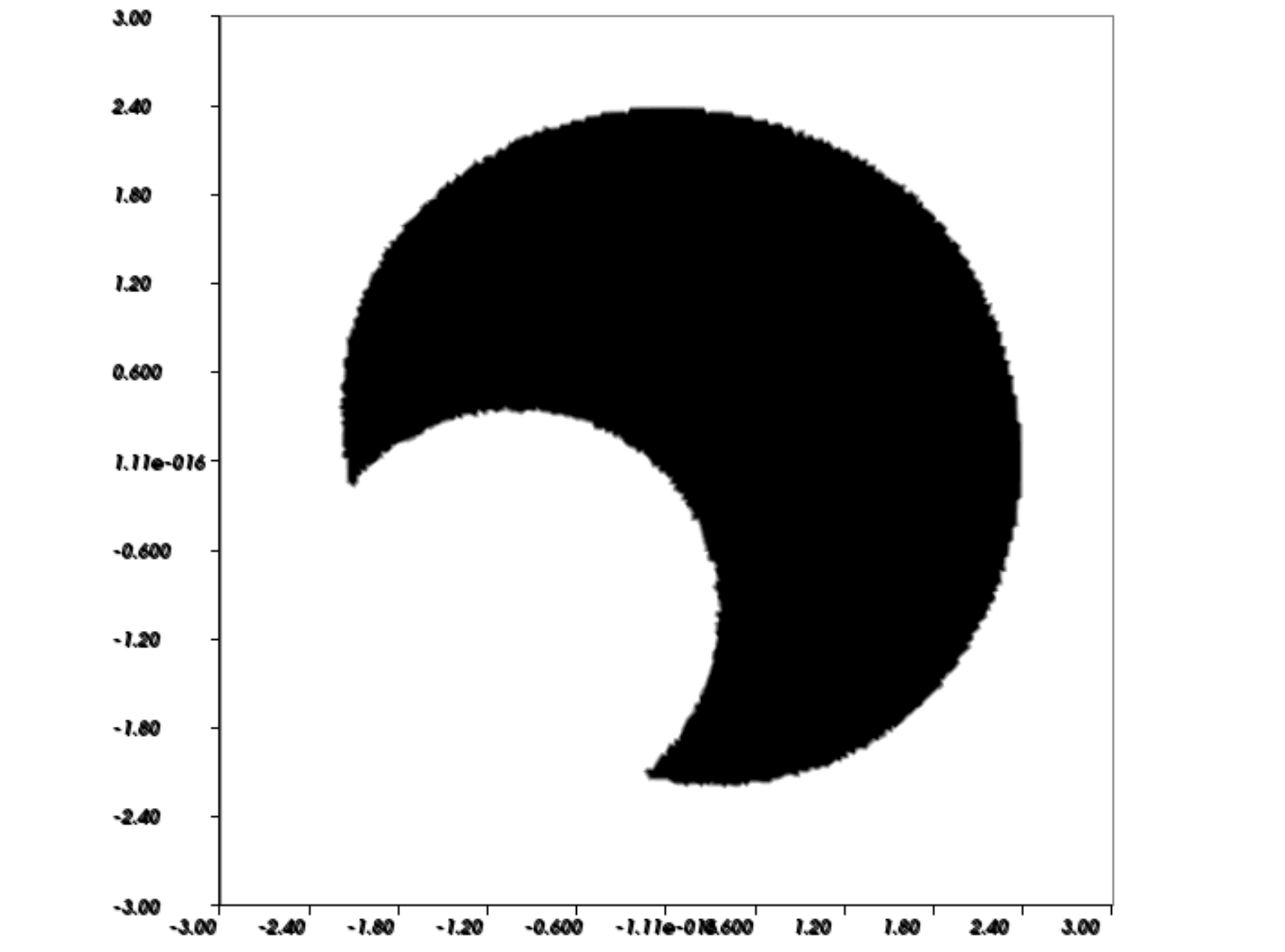}
\
\includegraphics[width=4.5cm]{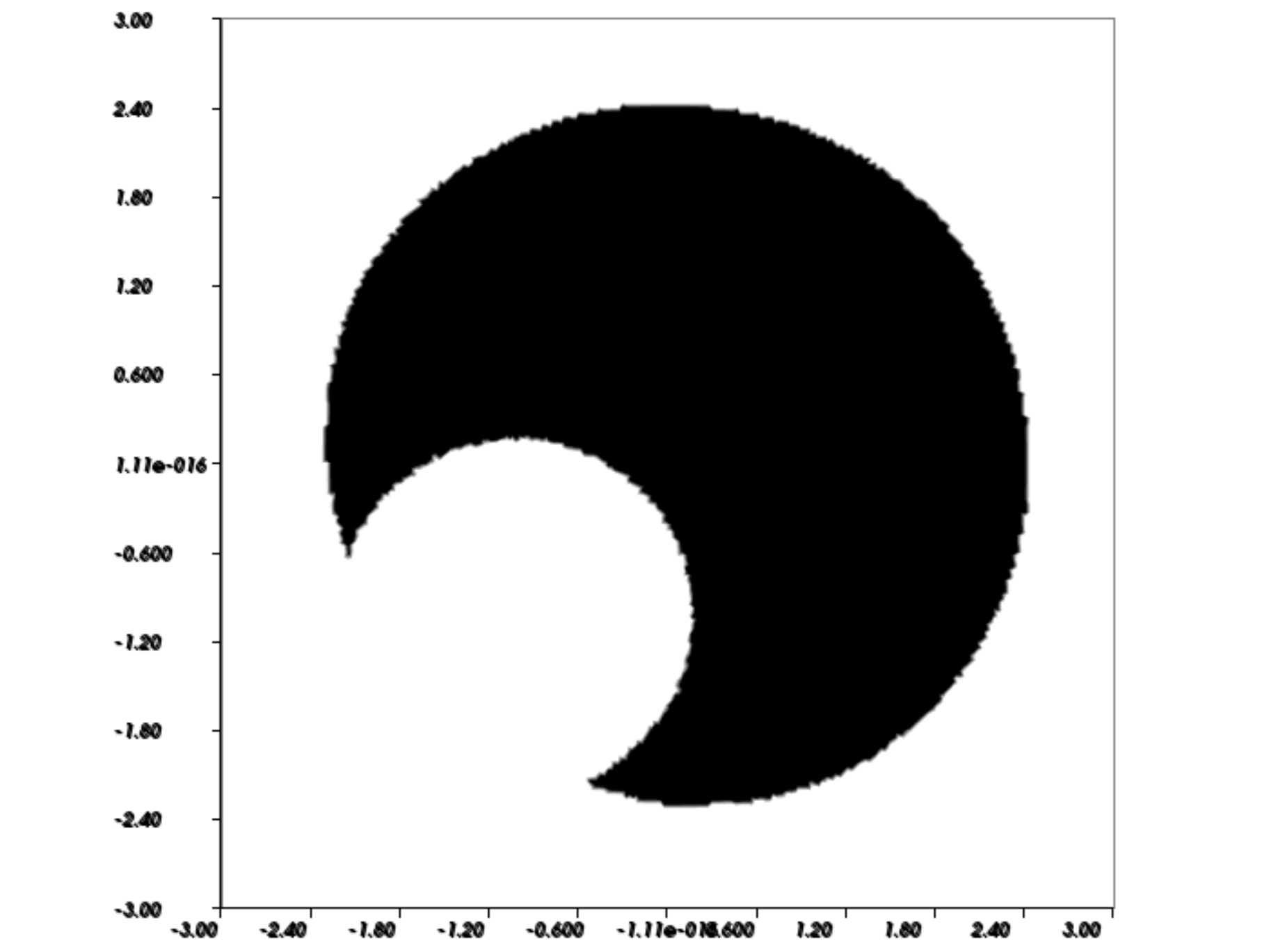}
\end{center}
\caption{Example 1. Initial domain $k=0$ (top, left), intermediate domains
during the line-search after $k=1$  
and the final domain $k=2$ (bottom, right).
\label{fig:ex1_xi}}
\end{figure}

We can observe in Figure \ref{fig:ex1_xi} the evolution of the
domain (both boundary and topological changes) and in
Table \ref{tab:ex1_J} the corresponding values of the
objective function. For $u_0=0$, we get $g_1=g_0$, but we have, for the cost functional, 
$\mathcal{J}_1 < \mathcal{J}_0$, since there is minimization with respect to the control $u$.
We do not plot in Figure \ref{fig:ex1_xi} the domain for $k=1$ because it
is the same as for 
$k=0$, but there is a column in Table \ref{tab:ex1_J} corresponding to $k=1$,
showing the evolution of the penalized cost.

\begin{table}[ht]
\begin{center}
\begin{tabular}{|c|c|c|c|c|c|c|c|}
\hline
iteration      & k=0    & k=1    &        &        &        &        & k=2    \\  \hline
$t_2$          & 220.87 & 171.13 & 155.60 & 149.47 & 129.19 &  67.60 &  90.50 \\  \hline
$t_3$          &  35.50 &  34.63 &  40.12 &  38.06 &  32.10 &  54.75 &  18.30 \\  \hline
$\mathcal{J}$  & 291.89 & 240.39 & 235.85 & 225.60 & 193.39 & 177.12 & 127.11 \\  \hline
\end{tabular}
\end{center}
\caption{Example 1. The computed objective function 
$\mathcal{J}=t_2+\frac{1}{\epsilon}t_3$.
The columns 4, 5, 6, 7 correspond to the
intermediate configurations
obtained during the line-search after $k=1$. The descent property is valid just for the total cost, on the last line.}
\label{tab:ex1_J}
\end{table}

\begin{table}[ht]
\begin{center}
\begin{tabular}{|c|c|c|c|c|c|c|}
\hline
iteration      & k=0   &       &       &       &       & k=2   \\  \hline
$t_2$          & 96.39 & 74.76 & 79.98 & 253.41& 46.59 & 56.62 \\  \hline
\end{tabular}
\end{center}
\caption{Example 1. The values of $t_2$ for the finite element solution
  of (\ref{2.1})-(\ref{2.2}) in the domains presented in Figure \ref{fig:ex1_xi}.}
\label{tab:ex1_t2}
\end{table}
For the solution of the elliptic problem (\ref{2.1})-(\ref{2.2})
in the computed domains $\Omega_g$, we obtain in fact the best value
$t_2=46.59$ (see Table \ref{tab:ex1_t2}), which is consistently better
than $t_2=67.60$ obtained for the
solution of (\ref{3.2})-(\ref{3.3}) in $D$, in the corresponding iteration
of the algorithm. This is due to the value of the penalization term $t_3$,
 which remians ``far'' from zero.
Such situations are frequent in penalization approaches for nonconvex
minimization problems.

\clearpage

\medskip
\noindent
\textbf{Example 2.}

We study now a case with $E\neq\emptyset$. The $D$, $y_d$, $f$, $\delta$
are the same as in Example 1. The observation domain $E$ is the disk of center
$(0,0)$ and radius $0.5$ and we take
$J(\mathbf{x})=\frac{1}{2}\left(y(\mathbf{x})-y_d(\mathbf{x})\right)^2$
and $j=0$. We fix $\epsilon=0.9$ and the other
numerical parameters are the same as in Example 1. Such a choice of a ``big'' 
penalization parameter (similar with the previous example) has the consequence
that the constraint (\ref{2.2})
is consistently relaxed and allows a large choice of descent
directions.

For $g_0(x_1,x_2)$, given by
$$
\max\left(
 (x_1+0.8)^2 + (x_2+0.8)^2 -1.8^2,
-(x_1+0.8)^2 - (x_2+0.8)^2 +0.6^2
\right) 
$$
we obtain as initial domain the ring of center $(-0.8,-0.8)$, exterior
radius $1.8$ and interior radius $0.6$.

In order to observe during the algorithm the restriction (\ref{2.5}), we use
the descent direction method
with projection, see \cite{Ciarlet2018}.
The descent direction is given by the Proposition \ref{prop:5.1}, case ii)
and the projection is computed as follows:
$\Pi(g)=g_E$ in $E$ and $\Pi(g)=g$ outside $E$,
where $g_E\in \mathcal{F}$ is such that $g_E(\mathbf{x})<0$ if and only if
$\mathbf{x}\in E$. In our test, $g_E(x_1,x_2)= (x_1)^2 + (x_2)^2 -0.5^2$.
The line search, with projection only for the parametrization function, is
$$
\lambda_k \in \arg\min_{\lambda >0}
\mathcal{J}\left( \Pi(G^k+\lambda R^k), U^k+\lambda V^k \right)
$$
and the next iteration is defined by
$$
G^{k+1}=\Pi(G^k+\lambda_k R^k),\quad U^{k+1}= U^k+\lambda_k V^k.
$$

\medskip

The initial guess for the control is $u_0=1$.

\begin{figure}[ht]
\begin{center}
\includegraphics[width=4.5cm]{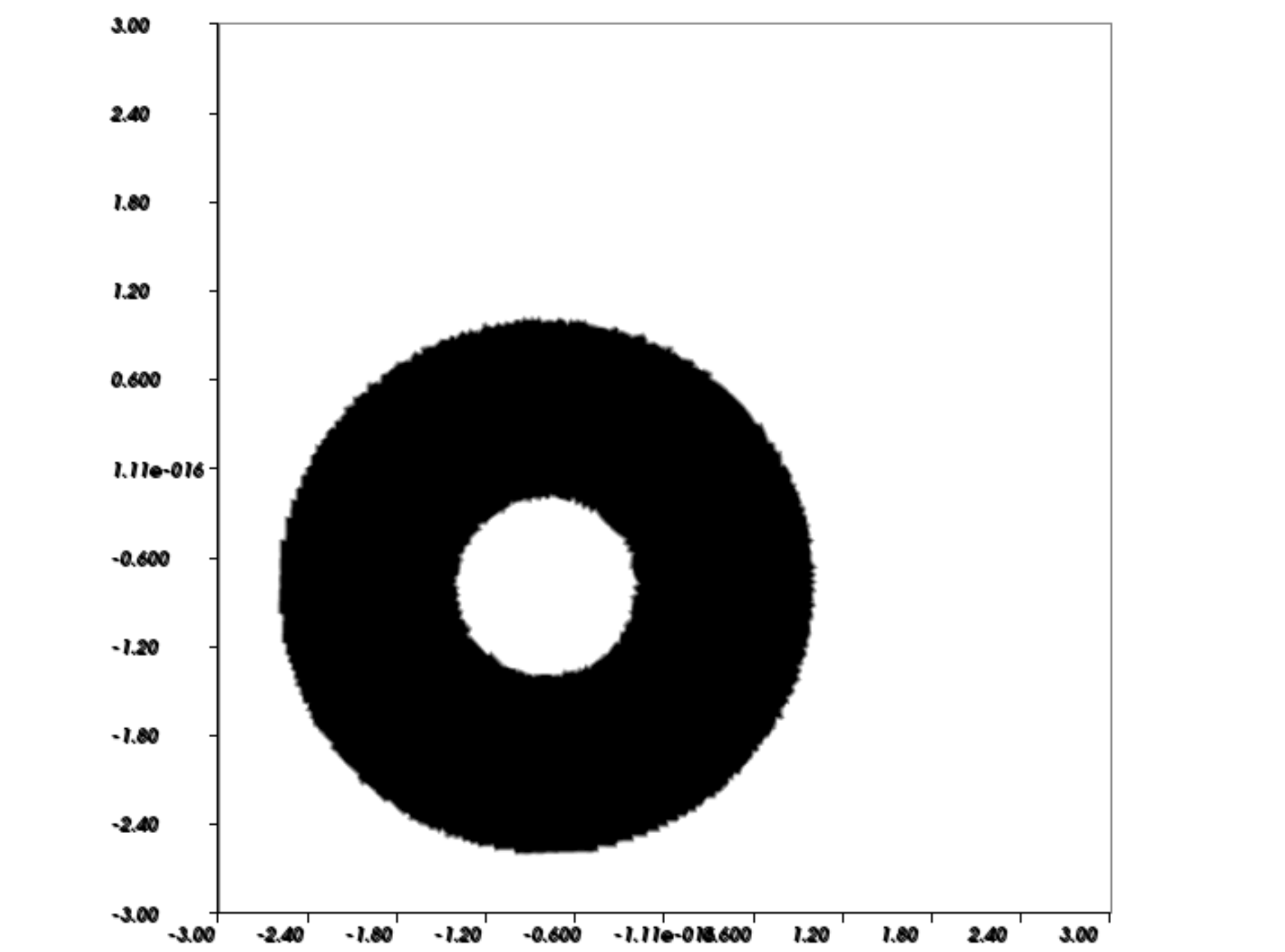}
\
\includegraphics[width=4.5cm]{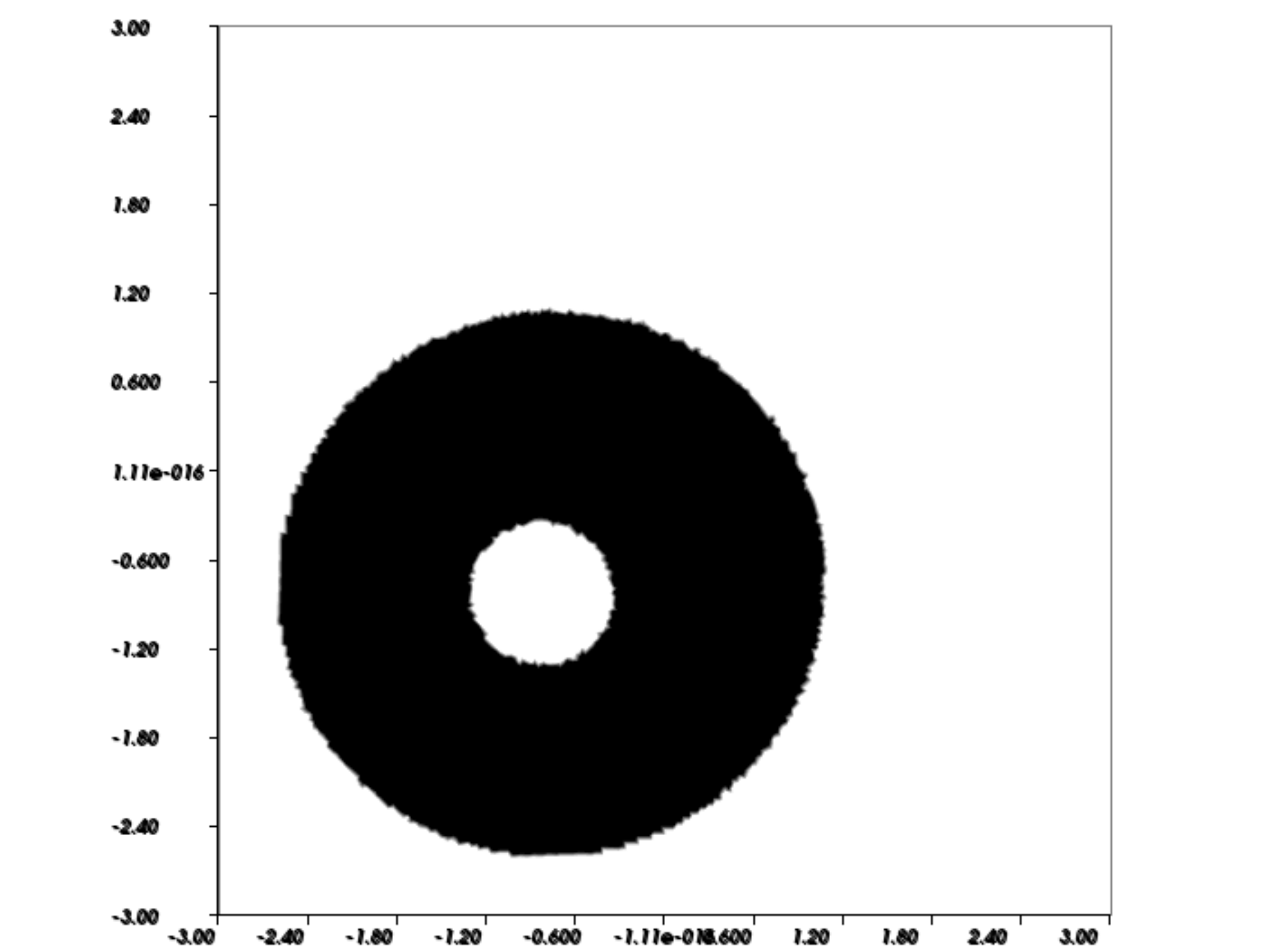}
\
\includegraphics[width=4.5cm]{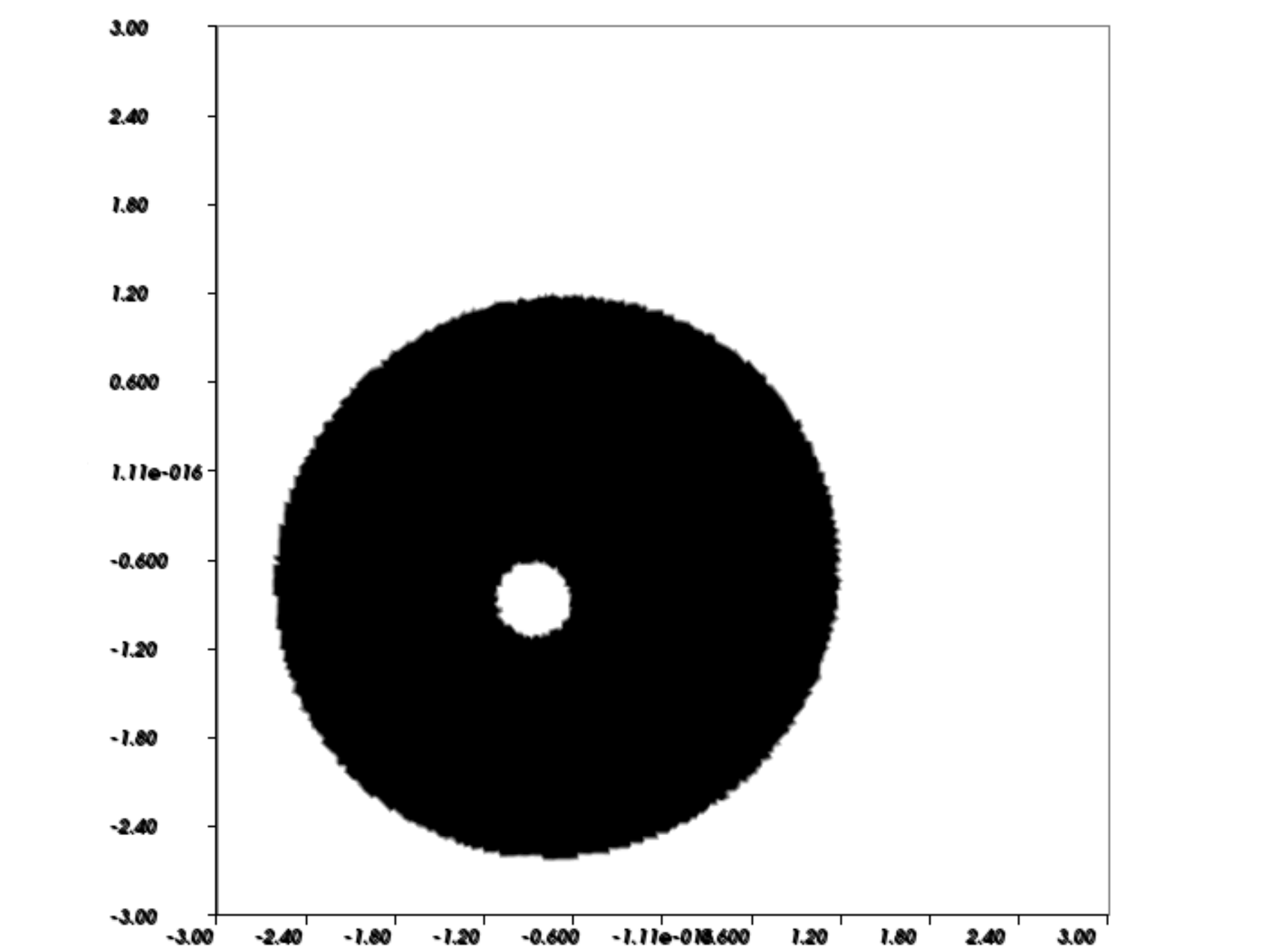}\\
\includegraphics[width=4.5cm]{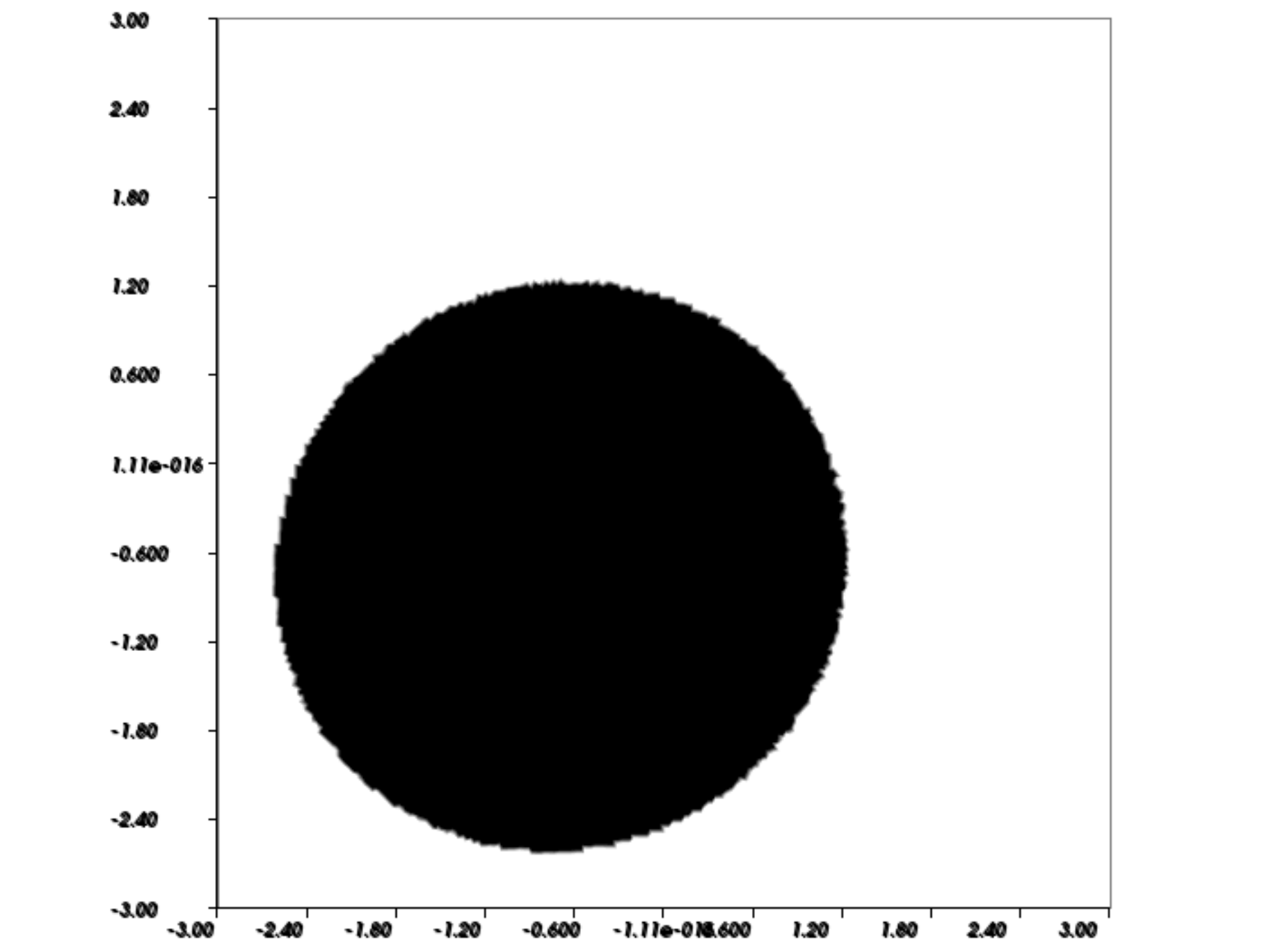}
\
\includegraphics[width=4.5cm]{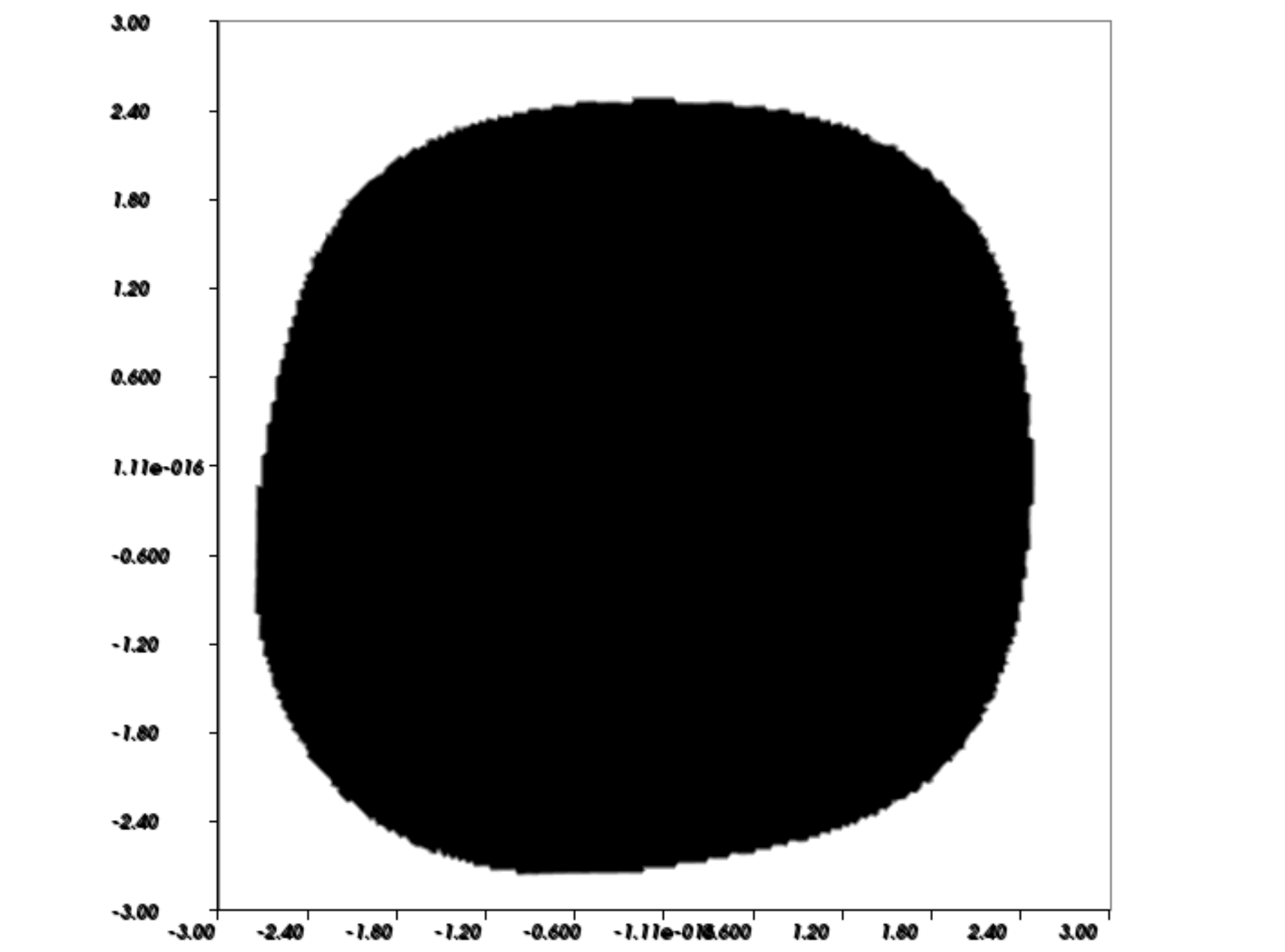}
\
\includegraphics[width=4.5cm]{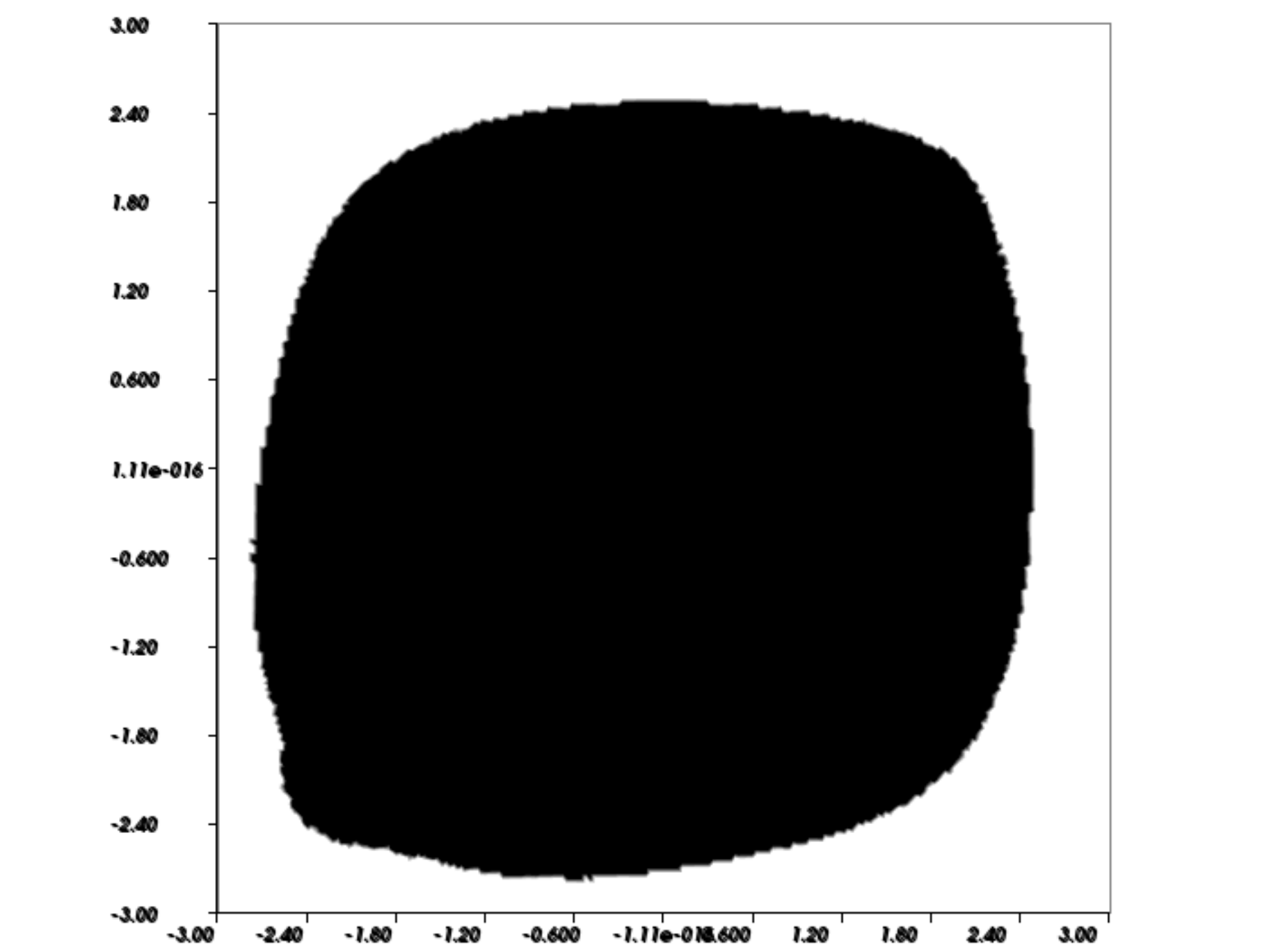}
\end{center}
\caption{Example 2. Domain for $k=0$ (top, left), intermediate domains
during the line-search after $k=0$,
domain for  $k=1$ (bottom, middle)
and the final domain for $k=2$ (bottom, right).
\label{fig:ex2_b_xi}}
\end{figure}

\begin{table}[ht]
\begin{center}
\begin{tabular}{|c|r|r|r|r|r|r|}
\hline
iteration      & k=0    &        &        &        & k=1    & k=2   \\  \hline
$t_1$          &   8.03 &   6.01 &   4.00 &   3.37 &   0.35 &  0.54 \\  \hline
$t_3$          & 234.91 & 218.34 & 204.47 & 198.08 & 193.56 & 57.42 \\  \hline
$\mathcal{J}$  & 269.05 & 248.62 & 231.20 & 223.46 & 215.42 & 64.35 \\  \hline
\end{tabular}
\end{center}
\caption{Example 2. The computed objective function 
$\mathcal{J}=t_1+\frac{1}{\epsilon}t_3$.
The columns 3, 4, 5 correspond to the
intermediate configurations
obtained during the line-search after $k=0$.}
\label{tab:ex2_b_J}
\end{table}

The domain evolution is presented in Figure \ref{fig:ex2_b_xi}
and the corresponding values of the
objective function are in Table \ref{tab:ex2_b_J}.

For the finite element solution of (\ref{2.1})-(\ref{2.2})
in the domains presented in Figure \ref{fig:ex2_b_xi}, we have reported $t_1$
in Table \ref{tab:ex2_b_t1}. Due to the low value of the initial cost, we
notice the oscillations around this value and the minimal cost is attained
already in the first step of the line search. The interpretation of the
penalization term is similar as in the previous example.
\begin{table}[ht]
\begin{center}
\begin{tabular}{|c|r|r|r|r|r|r|}
\hline
iteration      & k=0  &        &     &     & k=1 & k=2   \\  \hline
$t_1$          & 0.099& 0.00053& 0.11& 0.27& 0.51& 0.49  \\  \hline
\end{tabular}
\end{center}
\caption{Example 2. The values of $t_1$ for the finite element solution
  of (\ref{2.1})-(\ref{2.2}) in the domains presented in Figure \ref{fig:ex2_b_xi}.}
\label{tab:ex2_b_t1}
\end{table}


\clearpage



\begin{thebibliography}{99}
  
\bibitem{A}
G. Allaire, 
{Conception optimale de structures}, 
Volume 58 of Math\'ematiques \& Applications [Mathematics \& Applications].
Springer-Verlag, Berlin, 2007.

\bibitem{alst}
V. Arnautu, H. Langmach, J. Sprekels, D. Tiba, 
On the approximation and optimization of plates, 
{Numer.Funct.Anal. Optim.} 21 (2000) no.3-4, 337--354.

\bibitem{bst}
M. Barboteu, M. Sofonea, D. Tiba, 
The control variational method for beams in contact with deformable obstacles, 
{Z. Angew. Math. Mech.} 92, (2012) no.1, 25 – 40.

\bibitem{Bucur2005}
D. Bucur, G. Buttazzo,
{Variational methods in shape optimization problems}, 
Progress in Nonlinear Differential Equations and their Applications, 65. 
Birkhauser Boston, Inc., Boston, MA, 2005. 

\bibitem{Ciarlet2002}
P. G. Ciarlet,  
{The finite element method for elliptic problems}. 
Classics in Applied Mathematics, 40. Society for Industrial
and Applied Mathematics (SIAM), 
Philadelphia, PA, 2002.

\bibitem{Ciarlet2018}
P. G. Ciarlet, 
Introduction to Numerical Linear Algebra and Optimisation,
Cambridge University Press,
2018.

\bibitem{Delfour2001} 
M.C. Delfour, J.P. Zolesio, 
{Shapes and Geometries, Analysis, Differential Calculus and Optimization}, 
SIAM, Philadelphia, 2001.

\bibitem{Grisvard1985}
P. Grisvard, 
{Elliptic Problems in Nonsmooth Domains}. 
London, Pitman, 1985.

\bibitem{HMT2016} 
A. Halanay, C.M. Murea, D. Tiba, 
Existence of a steady flow of Stokes fluid past a linear
elastic structure using fictitious domain, 
{J. Math. Fluid Mech.}, 18 (2016), no. 2, 397--413.

\bibitem{HMT2018}
A. Halanay, C.M. Murea, D. Tiba, 
Extension theorems related to a fluid-structure interaction problem, 
{Bull. Math. Soc. Sci. Math. Roumanie}, 61 (2018) 417--437.

\bibitem{freefem++}
F. Hecht,
New development in FreeFem++. 
{J. Numer. Math.} {20} (2012) 251--265.  
\texttt{http://www.freefem.org}

\bibitem{Haslinger1996} 
J. Haslinger, P. Neittaanm\"aki, 
{Finite element approximation of optimal shape design}, 
J. Wiley \& Sons, New York,
1996.

\bibitem{Henrot2005}
A. Henrot, M. Pierre, {Variations et optimisation de formes.
Une analyse g\'eo\-m\'e\-trique},
Springer, 2005.

\bibitem{Hirsch2014}
M.W. Hirsch,  S. Smale, L.R. Devaney,
{Differential Equations, Dynamical Systems and an Introduction to Chaos},
Elsevier, Academic Press, San Diego (2014).

\bibitem{MAJ}
A. Maury, G. Allaire, and F. Jouve, 
Shape optimization with the level set method for
contact problems in linearised elasticity, 
{SMAI-Journal of computational mathematics}
3 (2017), pp. 249–292.

\bibitem{MT2019}
C.M. Murea, D. Tiba, 
Topological optimization via cost penalization, 
{Topological Methods in Nonlinear Analysis} 
Volume 54, No. 2B, (2019), 1023--1050.

\bibitem{MT2019a}
C.M. Murea, D. Tiba,
Optimization of a plate with holes,
{Computers and Mathematics with Applications}
77 (2019) 3010--3020.

\bibitem{N_Tiba2012} 
P. Neittaanm\"aki, D. Tiba, 
Fixed domain approaches in shape optimization problems, 
{Inverse Problems} 28 (2012) 1--35.

\bibitem{NP_Tiba2009} 
P. Neittaanm\"aki, A. Pennanen, D. Tiba, 
Fixed domain approaches in shape optimization problems with Dirichlet 
boundary conditions,  
{Inverse Problems} 25 (2009) 1--18.

\bibitem{NS_Tiba2006} 
P. Neittaanm\"aki, J. Sprekels, D. Tiba, 
{Optimization of elliptic systems. Theory and applications}, Springer, 
New York, 2006.

\bibitem{N_Tiba2015}
M.R. Nicolai, D. Tiba, 
Implicit functions and parametrizations in dimension three:
generalized solutions.
{Discrete Contin. Dyn. Syst.} 35 (2015), no. 6, 2701--2710. 

\bibitem{OF}
S. Osher and R. Fedkiw, 
{Level set methods and dynamic implicit surfaces}, 
Volume 153 of Applied
Mathematical Sciences. Springer-Verlag, New York, 2003.

\bibitem{OS}
S. Osher and J.A. Sethian,
Fronts propagating with curvature-dependent speed: algorithms based
on Hamilton-Jacobi formulations. 
{J. Comput. Phys.} 79 (1988), no. 1, 12--49.

\bibitem{Pironneau1984} 
O. Pironneau, 
{Optimal shape design for elliptic systems}, Springer, Berlin, 1984.

\bibitem{Pon}
L.S. Pontryagin, 
{Equations Differentielles Ordinaires}, MIR, Moscow, 1968

\bibitem{Raviart2004}
P.-A. Raviart and J.-M. Thomas, 
{Introduction \`a l'analyse num\'erique des \'equations aux d\'eriv\'ees partielles.} 
Dunod, 2004.


\bibitem{Sokolowski1992}
J. Sokolowski, J.P. Zolesio, 
{Introduction to Shape Optimization. Shape Sensitivity Analysis}, 
Springer, Berlin,
1992.


\bibitem{Tiba2013}
D. Tiba, 
The implicit function theorem and implicit parametrizations.
{Ann. Acad. Rom. Sci. Ser. Math. Appl.} 5 (2013), no. 1--2, 193--208. 

\bibitem{Tiba2018}
D. Tiba,  
Iterated Hamiltonian type systems and applications. 
{J. Differential Equations} 264 (2018), no. 8, 5465--5479. 

\bibitem{Tiba2018a}
D. Tiba, 
A penalization approach in shape optimization,
{Atti della Accademia  Peloritana dei
  Pericolanti - Classe di Scienze Fisiche, Matematiche e Naturali}
96 (2018), no. 1, A8.

\bibitem{Tiba2020}
D. Tiba,
Implicit parametrizations and applications in optimization and control,
{Mathematical Control and Related Fields} (2020), First online.
  
\end{thebibliography}
\end{document}